\documentclass[11pt]{article}
\usepackage{amsfonts}

\usepackage{amsmath}
\usepackage{graphicx}
\usepackage{makecell}
\usepackage{fancyhdr}

\newcommand{\eqnb}{\begin{equation}}
\newcommand{\eqne}{\end{equation}}

\newtheorem{The}{Theorem}

\newtheorem{Rem}{Remark}

\begin{document}

\title{A Mean-Field Matrix-Analytic Method for Bike Sharing Systems under
Markovian Environment}
\author{Quan-Lin Li and Rui-Na Fan \\
School of Economics and Management Sciences \\
Yanshan University, Qinhuangdao 066004, P.R. China}
\maketitle

\begin{abstract}
This paper proposes a novel mean-field matrix-analytic method in the study
of bike sharing systems, where the Markovian environment is constructed to
well express time-inhomogeneity and asymmetry of the processes that the
customers rent and return the bikes. To achieve effective computability of
this mean-field method, this paper provides a unified framework through
three basic steps: The first one is to deal with a major challenge
encountered in setting up the mean-field block-structured equations in more
general bike sharing systems. Here we provide an effective technique to
establish a necessary reference system which is a time-inhomogeneous queue
with block structure. The second step is to prove the asymptotic
independence (or \textit{propagation of chaos}) in terms of the martingale
limits. Note that the asymptotic independence ensures and supports that we
can construct a nonlinear QBD process such that the stationary probability
of problematic stations can be computed under a unified nonlinear QBD
framework. Finally, in the third step we use some numerical examples to show
effectiveness and computability of this mean-field matrix-analytic method,
and also give valuable observation on influence of some key parameters on
system performance. We hope the methodology and results given in this paper
are applicable in the study of more general large-scale bike sharing systems.

\vskip                                          0.5cm

\end{abstract}

\section{ Introduction}

In this paper, we apply the mean-field theory, combining Markov processes
with time-inhomogeneous queues, to study a complicated bike sharing system
with user's finite waiting rooms under a Markovian environment. To this end,
we develop a novel mean-field matrix-analytic method in the study of more
general bike sharing systems as follows: Setting up a block-structured
system of mean-field equations by constructing a reference system: The
time-inhomogeneous MAP$\left( t\right) $/MAP$\left( t\right) $/1/K+2L+1
queue; proving the asymptotic independence by means of the martingale
limits; establishing a nonlinear QBD process such that the fixed point can
be computed numerically; and using some numerical examples to give valuable
observation of influence of some key parameters on system performance.

During the last decades the bike sharing systems have emerged as a public
transport mode devoted to short trips, and they have widely been deployed in
more than 900 major cities around the world. So far the bike sharing systems
have been regarded as a promising solution to reduce traffic congestion,
parking difficulties, automobile exhaust pollution, transportation noise and
so on. For the history and survey papers, readers may refer to, such as,
DeMaio \cite{DeM:2009}, Shaheen \emph{et al}. \cite{Sha:2010}, Meddin and
DeMaio \cite{Med:2012}, Shu \emph{et al}. \cite{Shu:2013}, Fishman \emph{et
al}. \cite{Fis:2013}, Labadi \emph{et al}. \cite{Lab:2015}, Kaspi \emph{et al%
}. \cite{Kas:2017} and the references therein.

For design and operations of the bike sharing systems, it has become a basic
and interesting topic to assess and ensure the quality of service
(abbreviated as QoS) from a user's perspective, e.g., see Kaspi \emph{et al}%
. \cite{Kas:2017} for necessary interpretation. In general, the QoS of a
bike sharing system may be evaluated from two basic points: (a)\textit{\ The
bike non-empty}. Some bikes have been parked at the stations such that any
arriving customer can rent a bike from his entering station. (b)\textit{\
The parking non-full}. Some parking places (or lockers) become empty and
available so that a rider can immediately return his bike at a destination
station. Based on the two points, the bike-empty or parking-full stations
are called problematic stations, while the probability of problematic
stations can be used to measure the QoS of the bike sharing system. In
general, computing the probability of problematic stations is always very
difficult and challenging. To indicate the major reason of such
difficulties, from a physical point of view, Li \emph{et al}. \cite%
{LiFM:2016, LiF:2017} transformed the more general bike sharing system into
a complicated closed queueing network whose customers and nodes are all
virtual from the bikes, from the stations and from the roads, and provided
an effective method to compute the stationary probability of problematic
stations through deriving the product-form solution of joint stationary
distribution of queue lengths. On the other hand, it is worthwhile to note
that recent interesting research of bike sharing systems is also related to
the probability of problematic stations. Readers may refer to, for example,
the inventory management by Raviv and Kolka \cite{Rav:2013} and Schuijbroek
\emph{et al}. \cite{Sch:2013}, optimization of the bike fleet sizes by
Fricker and Gast \cite{Fri:2016}, and influence of the unusable bikes by
Kaspi \emph{et al}. \cite{Kas:2017}.

Little work has been done on numerical computation of the probability of
problematic stations through applications of queueing theory and Markov
processes. Important examples in the recent literature are classified as the
following two aspects. (a)\textit{\ Simple queues:} Leurent \cite{Leu:2012}
used the M/M/1/C queue to study a vehicle-sharing system. Schuijbroek \emph{%
et al}. \cite{Sch:2013} first evaluated the service level by means of the
transient distribution of the M/M/1/C queue, and then established some
optimal problems with respect to the vehicle routing. Raviv \emph{et al}. %
\cite{Rav:2012} and Raviv and Kolka \cite{Rav:2013} employed the transient
distribution of the M$\left( t\right) $/M$\left( t\right) $/1/C queue to
compute the expected number of bike shortages at each station. (b)\textit{\
Closed queueing networks:} Adelman \cite{Ade:2007} used a closed queueing
network to set up an internal pricing mechanism for managing a fleet of
service units, and provided a price-based policy for the vehicle
redistribution. George and Xia \cite{Geo:2011} used some simple closed
queueing networks to study the vehicle rental systems, and determined the
optimal number of parking places in each rental location. Waserhole \emph{et
al}. \cite{Was:2013} and Waserhole and Jost \cite{Was:2016} used a simple
closed queueing network, combining with the fluid approximation, to
establish the Markov decision models in order to determine the optimal
policy of the bike sharing system. Li \emph{et al}. \cite{LiFM:2016}
provided a unified framework of applying the closed queueing networks whose
product solution can be used to compute the stationary probability of
problematic stations. Li \emph{et al}. \cite{LiF:2017} further extended the
modeling method of \cite{LiFM:2016} to be able to study a more general bike
sharing system from two key factors: Markovian arrival processes and an
irreducible road graph. Under the heavy-traffic conditions, Li \emph{et al}. %
\cite{LiQF:2017} gave fluid and diffusion limits of the bike sharing system
with renewal user arrivals and general riding-bike times by means of
analyzing a multiclass closed queueing network.

If a bike sharing system contains $N$ stations and at most $N\left(
N-1\right) /2$ roads (or an irreducible road graph), then it can be
described as a virtual closed queueing network whose analysis is complicated
and difficult due to multiple virtual nodes (station or road) and many
parameters in this system. See Li \emph{et al}. \cite{LiF:2017, LiFM:2016}
for detailed interpretation. In this case, the mean-field theory should be
one of the best approximate methods for understanding dynamic behavior of
more general bike sharing systems. Note that the mean-field theory is a
popular approximate method in the study of complex physical systems, while
its applications contain at least two basic points: Focusing on a tagged
node, and computing mean-field parameters of this tagged node through weak
interactions among all the nodes. So far the available results of applying
the mean-field theory to the bike sharing systems have been very limited.
Fricker \emph{et al}. \cite{Fri:2012} first made a pioneering seminar work
for applying the mean-field theory to some heterogeneous bike sharing
systems. Since then, subsequent papers have been published on this theme.
Fricker and Gast \cite{Fri:2016} gave some simple mean-field models to study
a space-homogeneous bike sharing system in terms of the M/M/1/K queue, and
derived the closed-form solution both for the minimal proportion of
problematic stations and for the optimal fleet sizes. Fricker and Servel %
\cite{Fri:2016b} applied the mean-field theory to consider two-choice
regulation in heterogeneous closed networks, and then dealt with a bike
sharing system with multiple clusters. Fricker and Tibi \cite{Fri:2017}
first studied the central limit and local limit theorems for the independent
(non-identically distributed) random variables, which provide support on
analysis of a generalized Jackson network with product-form solution. Then
they used the limit theorems to give a stationary asymptotic analysis for
the locally space-homogeneous bike sharing systems. Li \emph{et al}. \cite%
{Li:2016b} applied the mean-field theory to discuss the bike sharing system
with more random factors through a time-inhomogeneous queue and a nonlinear
birth-death process, and numerically computed the fixed point which gives
performance analysis of the bike sharing system.

The purpose of this paper is to improve the mean-field theory to be able to
study more general bike sharing systems from two key factors: (a) A
Markovian environment is constructed to well express time-inhomogeneity and
asymmetry of the processes that the customers rent and return the bikes; and
(b) user's finite waiting rooms are added to the stations such that the
probability of problematic stations can be reduced greatly. From
mathematical modeling and analysis, both the Markovian environment and the
user's finite waiting rooms make analysis of the bike sharing systems more
difficult and challenging. In addition, it is worthwhile to note that
introduction of the Markovian environment motivates us to improve the
mean-field theory to be able to set up the mean-field block-structured
equations, to prove the asymptotic independence with block structure, and to
establish a nonlinear Markov process which is used to compute the fixed
point. Based on this, we develop a mean-field matrix-analytic method in the
study of bike sharing systems. As a nearby research of this paper, Li and
Lui \cite{Li:2016c} applied the mean-field theory to discuss a
block-structured supermarket model. Li \cite{Li:2016a} provided a unified
block-structured framework for the mean-field theory of stochastic big
networks with weak interactions.

For the mean-field theory of stochastic networks, readers may refer to, such
as, two books by Liggett \cite{Lig:1985} and Chen \cite{Chen:2004}, two
survey papers by Sznitman \cite{Szn:1989} and Benaim and Le Boudec \cite%
{Ben:2008}. Since the mean-field theory was first applied to the study of
large-scale parallel queues (for example, supermarket models and work
stealing models) by Vvedenskaya \emph{et al}. \cite{Vve:1996} and
Mitzenmacher \cite{Mit:1996}, subsequent papers have been published on this
theme, among which see Turner \cite{Tur:1998}, Martin and Suhov \cite%
{Mar:1999}, Graham \cite{Gra:2000, Gra:2005}, Gast and Gaujal \cite%
{Gas:2010, Gas:2011}, Li \emph{et al}. \cite{Li:2014, Li:2015}, Li and Lui %
\cite{Li:2016c}, Li \cite{Li:2016a}, Mukhopadhyay \cite{Muk:2017} and the
references therein. On the other hand, the QBD processes often provide a
useful mathematical tool for studying stochastic models such as queueing
systems, manufacturing systems, communication networks and healthcare
systems. Readers may refer to Chapter 3 of Neuts \cite{Neu:1981}, Latouche
and Ramaswami \cite{Lat:1999}, Li \cite{Li:2010} and references therein.

The main contributions of this paper are threefold. The first one is to
propose a novel mean-field matrix-analytic method in the study of bike
sharing systems. Note that this new method can effectively improve the
descriptive and computational ability of the mean-field theory under a
unified framework of nonlinear Markov processes, e.g., see Li \cite{Li:2016a}
for a detailed discussion. To demonstrate such an ability by using examples,
in a bike sharing system we first introduce two key factors: The Markovian
environment and the user's finite waiting rooms. Then we show that the two
factors may result in some major challenges when applying the mean-field
theory to the bike sharing system, for example, it is always very difficult
to set up a block-structured system of mean-field equations due to existence
of the Markovian environment. To overcome difficulty of the block structure,
Subsection 3.1 provides an effective technique for establishing a necessary
reference system: A time-inhomogeneous queue MAP$\left( t\right) $/MAP$%
\left( t\right) $/1/K+2L+1. At the same time, the other key points of
applying the mean-field matrix-analytic method are also discussed as
follows: ($i$) Section 4 proves the asymptotic independence by means of the
martingale limits, ($ii$) Section 5 establishes a nonlinear QBD process such
that the fixed point can be computed numerically, and ($iii$) Section 6 uses
the fixed point to evaluate performance measures of the bike sharing system,
and specifically, to compute the stationary probability of problematic
stations.

The second contribution of this paper is to introduce the Markovian
environment, which can well express time-inhomogeneity and asymmetry of the
processes that the customers rent and return the bikes in a bike sharing
system. To our best knowledge, it is the first time that a Markovian
environment is constructed in the end of Section 2 by means of useful
information which arises from the rate volatility of the processes that the
customers rent and return the bikes within one period (i.e., one day), where
a fluctuating law of three peaks in the bike-rented (or bike-returned)
processes is refined from the practical data of the tagged station of the
bike sharing system, e.g. see Figure 1. On the other hand, the user's finite
waiting rooms are added into some stations, and they enhance flexibility and
ability of the bike sharing system such that the probability of problematic
stations is reduced greatly. Thus the QoS of the bike sharing system can be
promoted effectively by means of adding the user's finite waiting rooms,
e.g., see Figures 9 and 10 for some numerical analysis. Note that the user's
finite waiting rooms were first introduced and discussed by Leurent \cite%
{Leu:2012} through the M/M/1/C+K queue in which only one isolated station is
observed and analyzed, and the results obtained from the isolated station
were used to provide a coarse-grained approximation for performance
evaluation of the bike sharing system. Differently from Leurent \cite%
{Leu:2012}, this paper analyzes a total network of the bike sharing system
by means of the mean-field theory, where the nodes with finite waiting rooms
may have a variety of weak interactions. The third contribution is to use
some numerical examples to show effectiveness and computability of this
mean-field matrix-analytic method, and to show how some key parameters
influence performance measures of this bike sharing system. Therefore, we
gain new insights on understanding nonlinear dynamics, inhomogeneous nature
and interesting performance of the bike sharing systems, and hope the
methodology and results given in this paper are applicable in the study of
more general large-scale bike sharing systems.

The remainder of this paper is organized as follows. In Section 2, we
describe a large-scale bike sharing system with $N$ identical stations and
with user's finite waiting rooms under Markovian environment. Furthermore,
we provide a method to construct a Markovian environment by means of a
fluctuating rate law of three peaks within a period, which well expresses
time-inhomogeneity and asymmetry of the processes that the customers rent
and return the bikes. In Section 3, we first introduce an empirical measure
process to express the states of this bike sharing system. Then we use a
probability-analytic method to establish a necessary reference system: A
time-inhomogeneous queue MAP$\left( t\right) $/MAP$\left( t\right) $%
/1/K+2L+1 by means of the mean-field theory. This help us to set up a
block-structured system of mean-field equations. In Section 4, we apply the
martingale limit theory to prove the asymptotic independence (or propagation
of chaos) of the bike sharing system. In Section 5, we discuss the fixed
point of the block-structured system of limiting mean-field equations, and
provide a nonlinear QBD process to compute the fixed point. Furthermore, we
study the limiting interchangeability as $N\rightarrow \infty $ and $%
t\rightarrow +\infty $. In Section 6, we give six numerical examples to
investigate the performance measures, and show how some key parameters
influence system performance. Some concluding remarks are given in Section 7.

\section{Model Description}

In this section, we describe a large-scale bike sharing system with $N$
identical stations and with user's finite waiting rooms under Markovian
environment, and list operations mechanism, system parameters, model
notation and necessary interpretation. Furthermore, we give a detailed
discussion on how to construct a Markovian environment by means of useful
information arose from the rate volatility of the process that the customers
rent and return the bikes within one period (i.e., one day), where a
fluctuating law of three peaks in the bike-rented (or bike-returned)
processes is refined according to practical dynamics of the bike sharing
systems.

In a bike sharing system, a customer first arrives at a station, rents a
bike, and uses it for a while; then he returns the bike to a destination
station. Once the customer finishes his trip and returns the bike to a
station, he immediately leaves the bike sharing system. Based on this, Li
\emph{et al}. \cite{LiFM:2016} first described a practical bike sharing
system as a complicated closed queueing network with virtual customers
(bikes) and two classes of virtual nodes (stations and roads). Since then,
subsequent papers have been published on this theme. Li \emph{et al}. \cite%
{LiF:2017} extended the model in \cite{LiFM:2016} from two key factors:
Markovian arrival processes and an irreducible road graph; while Li \emph{et
al}. \cite{LiQF:2017} gave fluid and diffusion limits of the bike sharing
system with renewal user arrivals and general riding-bike times.

Although it is seen from the closed queueing networks that the stationary
probability of problematic stations can formally be computed by means of the
product-form solution of the joint stationary distribution of queue lengths,
there still exist some calculable difficulties which arise from a rather
complicated expression of the routing matrix of the corresponding virtual
closed queueing network. Therefore, it is not only important for theoretic
investigations but also necessary for engineering applications to provide
some effective approximate techniques, for example, the mean-field theory of
the bike sharing systems. Note that the mean-field theory has a key
advantage that focuses on analyzing only one node with mean-field parameters
whose basic information is sourced from the weak interactions among all the
nodes of a network system. Thus performance measures of this node with
mean-field parameters can be obtained easily, and they are used to well
approximate that of the total network system. Thus the mean-field results
can be used to show how some key system parameters influence the stationary
probability of problematic stations so that the QoS of the bike sharing
system can be evaluated from such a mean-field approximation.

Based on the above analysis, this paper extends the mean-field theory to be
able to deal with more general bike sharing systems, and further provides a
novel mean-field matrix-analytic method in the study of bike sharing
systems. To this end, we describe a large-scale bike sharing system, and
list its operations mechanism, system parameters, model notation and
necessary interpretation as follows:

\textbf{(1) The }$N$\textbf{\ identical stations:} To use the mean-field
theory, we assume that the large-scale bike sharing system contains $N$
identical stations; and at the initial time $t=0$, every station has $C$
bikes and $K$ parking places, in which $1\leq C\leq K<\infty $. Every
station continuously operates through either renting a bike or returning a
bike, so the number of bikes in a tagged station can be regarded as a
queueing process.

\textbf{(2) Adding user's finite waiting rooms to the stations:} To decrease
the probability of problematic stations when customers are sufficient in
this system, it is an efficient technique to add a user's finite waiting
room at each station. The waiting room have $L$ waiting places, each of
which is occupied by only a customer when either he can not rent a bike from
the tagged station or he can not return his bike to the tagged station. In
general, each finite waiting room has two useful purposes: (a) When a
customer arrives at an empty station in which no bike can be rented, either
he enters a waiting place in order to wait for a future available bike with
probability $\alpha \in \left[ 0,1\right] $, or he immediately leaves the
bike sharing system with probability $1-\alpha $. (b) When a riding-bike
customer completes his trip and enters a full station in which no empty
parking place is available, either he enters a waiting place in order to
wait for a future available parking place with probability $\beta \in \left[
0,1\right] $, or he immediately rides his bike to another station in order
to return the bike with probability $1-\beta $. Here, we must explain that
the riding-bike customer must return his bike to any station, and then he
can leave the bike sharing system, because each bike is an indispensable
public equipment and cannot be lost or become a personal property.

\textbf{(3) The Markovian environment: }In this bike sharing system, the
arrival and travel processes are influenced (or controlled) by a Markovian
environment, which is a continuous-time irreducible positive-recurrent
Markov process whose infinitesimal generator of size $m$ is given by

\begin{equation*}
\mathbf{W}=\left(
\begin{array}{cccc}
w_{1,1} & w_{1,2} & \cdots & w_{1,m} \\
w_{2,1} & w_{2,2} & \cdots & w_{2,m} \\
\vdots & \vdots & \ddots & \vdots \\
w_{m,1} & w_{m,2} & \cdots & w_{m,m}%
\end{array}%
\right) ,
\end{equation*}%
where $w_{i,i}<0$ for $1\leq i\leq m$; $w_{i,j}\geq 0$ for $1\leq i,j\leq m$
and $i\neq j$; $\sum_{j=1}^{m}w_{i,j}=0$ for $1\leq i\leq m$. At the same
time, we denote by $\theta $ the stationary probability vector of the Markov
process $\mathbf{W}$, that is, $\theta \mathbf{W}=0$ and $\theta e=1$, where
$e$ is a column vector of ones. Based on the Markov process $\mathbf{W}$,
now we describe the arrival processes and the riding-bike times as follows:

\textbf{(3.1) The arrival processes: }If the Markovian environment is at
State $j$, then the arrivals of customers at the bike sharing system from
outside are a Poisson process with arrival rate $N\lambda _{j}$ for $1\leq
j\leq m$.

\textbf{(3.2) The riding-bike times:} If the Markovian environment is at
State $j$, then the riding-bike time that a customer rides a bike from one
station to another is exponential with travel rate $\mu _{j}$ for $1\leq
j\leq m$.

\textbf{(4) The leaving principle: }Once a customer finishes his trip and
returns his bike to any station, he immediately leaves the bike sharing
system.

We assume that all the random variables defined above are independent of
each other.

When observing a tagged station in the bike sharing system, the finite
waiting room and the Markovian environment play a key role in queueing
analysis of this tagged station. To explain this, the queueing structure of
the tagged station is depicted in Figure 1.

\begin{figure}[ptb]
\setlength{\abovecaptionskip}{0.cm} \setlength{\belowcaptionskip}{-0.cm} %
\centering              \includegraphics[width=10cm]{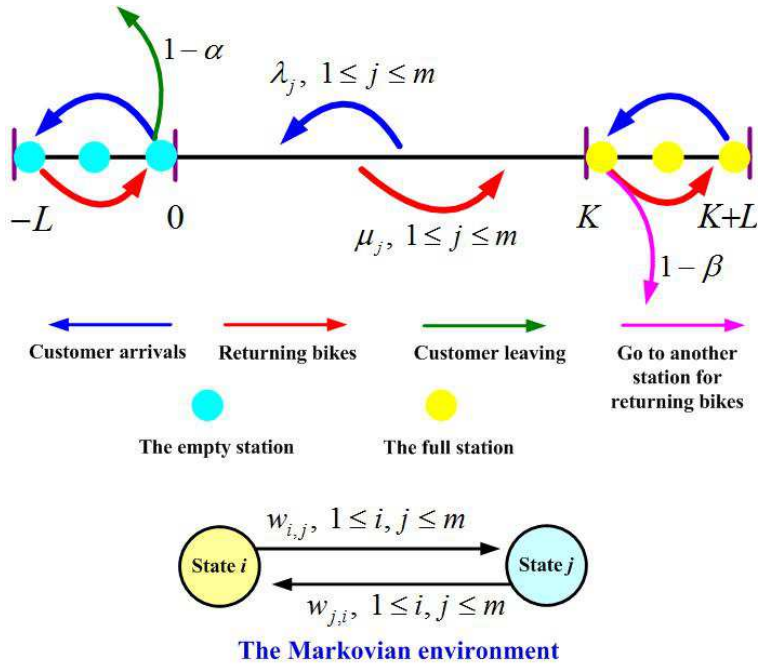} \newline
\caption{The queueing structure of any station in this bike-sharing system}
\label{figure:figure-1}
\end{figure}

\begin{Rem}
On the one hand, the assumption of the $N$ identical stations is used to
guarantee applicability of the mean-field theory (that is, the
multi-dimensional Markov process is exchangeable). On the other hand, from a
practical point of view, the stations in a major city are also designed as
almost the same, for example, Hangzhou has over 4000 stations, and each
station contains about 30 bikes.
\end{Rem}

\begin{Rem}
In some major cities, there are always many bikes and customers distributed
at various stations to support the short trips. To improve the quality of
service (or to decrease the probability of problematic stations), an
adscititious waiting room designed to add at each station of a bike sharing
system is always effective and useful. Leurent \cite{Leu:2012} first
proposed such an idea of adscititious waiting rooms, and discussed the
queueing process of only one isolated station by means of the M/M/1/C queue.
Differently from Leurent \cite{Leu:2012}, this paper applies the mean-field
theory to analyse such a network system of the bike sharing system with
user's waiting rooms, and then numerically compares performance measures of
the bike sharing systems with or without user's finite waiting rooms, e.g.,
see Figure 9 for more details.
\end{Rem}

In the remainder of this section, we give an interesting idea that
constructs a Markovian environment with seven states to be able to well
express the time-inhomogeneity and asymmetry of the processes that the
customers rent and return the bikes. To our best knowledge, it is the first
time that a Markovian environment is constructed by means of useful
information which arises from the rate volatility of the processes that the
customers rent and return the bikes within a day, where a fluctuating law of
three peaks in bike rented (or returned) rates within a day is depicted in
Figure 2.

\begin{figure}[ptb]
\setlength{\abovecaptionskip}{0.cm} \setlength{\belowcaptionskip}{-0.cm} %
\centering              \includegraphics[width=13cm]{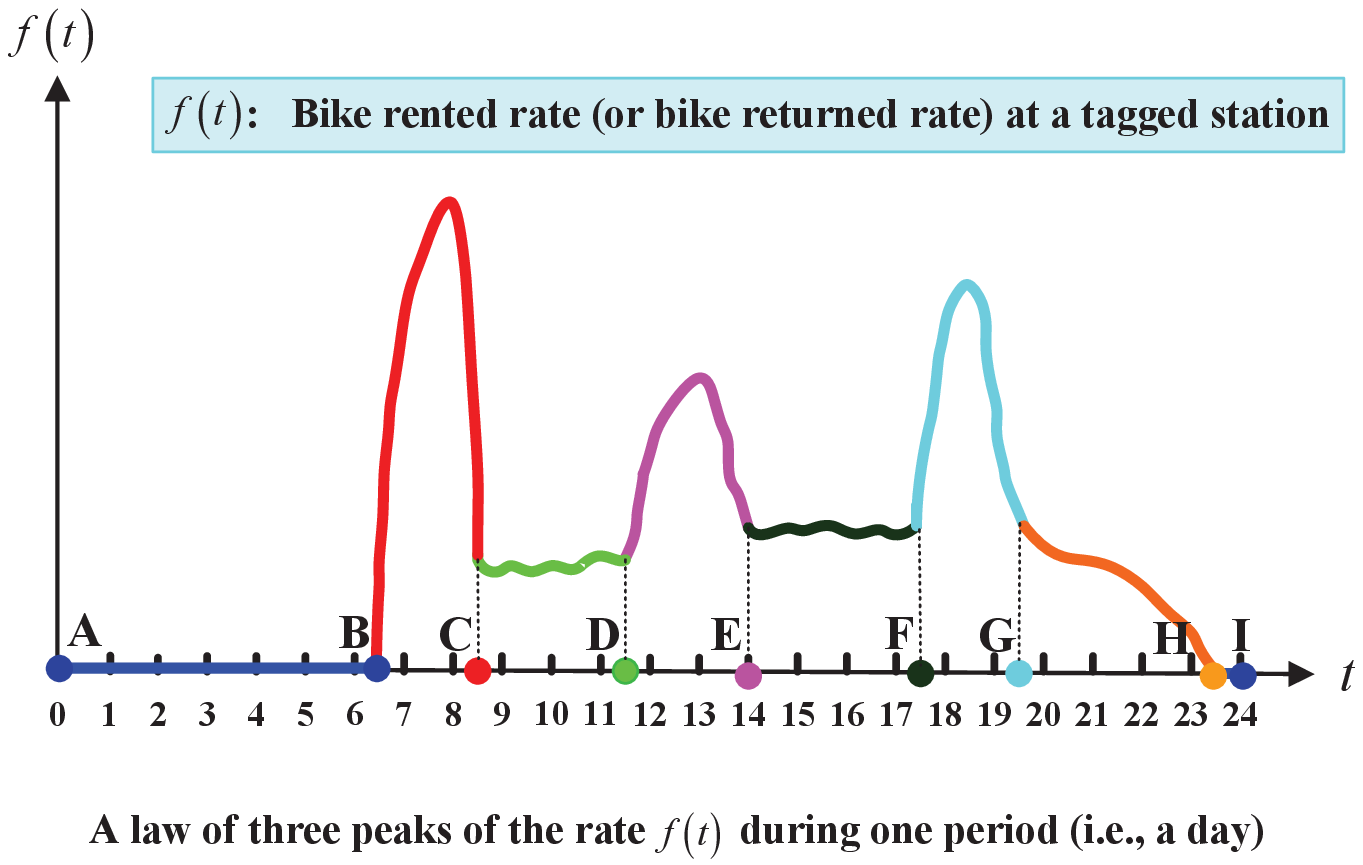} \newline
\caption{A fluctuating law of three peaks for bike rented (or returned)
rates within a day}
\label{figure:figure-2-1}
\end{figure}

Based on Figure 2, we segment 24 hours of one day into 7 parts as follows:%
\begin{align*}
\text{Part one}& =[0,6.5)\cup \lbrack 23.5,24),\text{ \ \ \ Part two}%
=[6.5,8.5), \\
\text{Part three}& =[8.5,11.5),\text{ \ \ \ \ \ \ \ \ \ \ \ \ \ \ \ Part four%
}=[11.5,14), \\
\text{Part five}& =[14,17.5),\text{ \ \ \ \ \ \ \ \ \ \ \ \ \ \ \ Part six}%
=[17.5,19.5), \\
\text{Part seven}& =[19.5,23.5).
\end{align*}%
Since the seven segmented parts within 24 hours demonstrate a stable regular
structure of periodical dynamics of the bike sharing systems, Part $i$ may
be regarded as State $i$ of a Markovian environment. Note that the Markovian
environment can be expressed by an irreducible Markov process whose state
transition relations are described as%
\begin{equation*}
\text{\textbf{State} }\mathbf{1}\rightarrow \text{State }2\rightarrow \text{%
State }3\rightarrow \text{State }4\rightarrow \text{State }5\rightarrow
\text{State }6\rightarrow \text{State }7\rightarrow \text{\textbf{State }}%
\mathbf{1}.
\end{equation*}%
Thus the infinitesimal generator of the Markovian environment is given by%
\begin{equation*}
\mathbf{W}=\left(
\begin{array}{ccccccc}
-x_{1} & x_{1} &  &  &  &  &  \\
& -x_{2} & x_{2} &  &  &  &  \\
&  & -x_{3} & x_{3} &  &  &  \\
&  &  & -x_{4} & x_{4} &  &  \\
&  &  &  & -x_{5} & x_{5} &  \\
&  &  &  &  & -x_{6} & x_{6} \\
x_{7} &  &  &  &  &  & -x_{7}%
\end{array}%
\right) .
\end{equation*}%
To compute the undetermined numbers $x_{i}$ for $1\leq i\leq 7$, we first
take the stationary probability vector $\theta =\left( \theta _{1},\theta
_{2},\ldots ,\theta _{7}\right) $ according to the time length ratios of the
seven segmented parts within 24 hours. Thus it is seen from Figure 2 that%
\begin{equation*}
\theta _{1}=\frac{7}{24},\theta _{2}=\frac{2}{24},\theta _{3}=\frac{3}{24}%
,\theta _{4}=\frac{2.5}{24},\theta _{5}=\frac{3.5}{24},\theta _{6}=\frac{2}{%
24},\theta _{7}=\frac{4}{24}.
\end{equation*}%
Let $x_{7}=1$. Then it follows from $\theta \mathbf{W}=0$ that%
\begin{equation*}
x_{i}=\frac{\theta _{7}}{\theta _{i}},\text{ \ \ }1\leq i\leq 7.
\end{equation*}

Now, we provide an average approximate method to determine the bike rented
and returned rates which are controlled by the states of the Markovian
environment. To this end, let $f_{\text{rent}}\left( t\right) $ and $f_{%
\text{return}}\left( t\right) $ be the instantaneous rates of the bike
rented and returned processes at time $t\geq 0$, respectively. Note that
each of the two functions $f_{\text{rent}}\left( t\right) $ and $f_{\text{%
return}}\left( t\right) $ is referred to the fluctuating law of three peaks
depicted in Figure 2, or they can be approximately evaluated by means of the
associated data collected from system operations. Based on this, we set up
\begin{align*}
\lambda _{1}& =0,\text{ }\mu _{1}=0; \\
\lambda _{2}& =\frac{1}{2}\int_{6.5}^{8.5}f_{\text{rent}}\left( t\right)
\text{d}t,\text{ }\mu _{2}=\frac{1}{2}\int_{6.5}^{8.5}f_{\text{return}%
}\left( t\right) \text{d}t; \\
\lambda _{3}& =\frac{1}{3}\int_{8.5}^{11.5}f_{\text{rent}}\left( t\right)
\text{d}t,\text{ }\mu _{3}=\frac{1}{3}\int_{8.5}^{11.5}f_{\text{return}%
}\left( t\right) \text{d}t; \\
\lambda _{4}& =\frac{1}{2.5}\int_{11.5}^{14}f_{\text{rent}}\left( t\right)
\text{d}t,\text{ }\mu _{4}=\frac{1}{2.5}\int_{11.5}^{14}f_{\text{return}%
}\left( t\right) \text{d}t; \\
\lambda _{5}& =\frac{1}{3.5}\int_{14}^{17.5}f_{\text{rent}}\left( t\right)
\text{d}t,\text{ }\mu _{5}=\frac{1}{3.5}\int_{14}^{17.5}f_{\text{return}%
}\left( t\right) \text{d}t; \\
\lambda _{6}& =\frac{1}{2}\int_{17.5}^{19.5}f_{\text{rent}}\left( t\right)
\text{d}t,\text{ }\mu _{6}=\frac{1}{2}\int_{17.5}^{19.5}f_{\text{return}%
}\left( t\right) \text{d}t; \\
\lambda _{7}& =\frac{1}{4}\int_{19.5}^{23.5}f_{\text{rent}}\left( t\right)
\text{d}t,\text{ }\mu _{7}=\frac{1}{4}\int_{19.5}^{23.5}f_{\text{return}%
}\left( t\right) \text{d}t.
\end{align*}%
In general, the two functions $f_{\text{rent}}\left( t\right) $ and $f_{%
\text{return}}\left( t\right) $ can approximately be given by means of some
statistical models to deal with the practical data in a tagged station of
the bike sharing system.

\begin{Rem}
By using the fluctuating law of three peaks depicted in Figure 2, we choose
seven different states to construct a Markovian environment, which is
related to real-time dynamics of the bike sharing system. Similarly, we may
also set up a fluctuating law of $n$ peaks to construct a Markovian
environment of $2n+1$ states. Note that the approximate accuracy of such
modeling can be improved effectively as the number $n$ increases. On the
other hand, it is necessary to mention that the Markovian environment
constructed from a fluctuating law of $n$ peaks will be very useful in the
study of stochastic periodic systems because the difficult periodic dynamic
system is transformed to an easy Markov system, for example, ride sharing
systems, healthcare systems, transportation networks, periodic retail
systems, wind power system and so forth.
\end{Rem}

\section{Mean-Field Equations}

In this section, we first introduce an empirical measure process to express
the states of this bike sharing system. Then we provide a
probability-analytic method, combining with the mean-field theory, to
establish a necessary reference system: A time-inhomogeneous queue MAP$%
\left( t\right) $/MAP$\left( t\right) $/1/K+2L+1. This help us to set up a
block-structured system of mean-field equations.

Now, we introduce an empirical measure process to express the states of this
bike sharing system.

Let $X_{i}^{\left( N\right) }\left( t\right) $ and $J\left( t\right) $ be
the number of customers in Station $i$ and the state of the Markovian
environment at time $t$, respectively. Then $\mathbf{X}%
=\{(X_{1}^{(N)}(t),X_{2}^{(N)}(t),\ldots ,X_{N}^{(N)}(t);$\newline
$J(t)):t\geq 0\}$ is an $Nm$-dimensional Markov process due to the
assumptions on the Poisson arrivals, the exponential travel times and the
Markovian environment. Note that analysis of such an $Nm$-dimensional Markov
process is always difficult due to the ``\textit{State Space Explosion}''.
Therefore, it is necessary to introduce an empirical measure process.

For the Markov process $\mathbf{X=}\left\{ \left( X_{1}^{\left( N\right)
}\left( t\right) ,X_{2}^{\left( N\right) }\left( t\right) ,\ldots
,X_{N}^{\left( N\right) }\left( t\right) ;J\left( t\right) \right)
:t\geq0\right\} $, the empirical measure is defined as
\begin{equation*}
Y_{k,j}^{\left( N\right) }\left( t\right) =\frac{1}{N}\sum\limits_{n=1}^{N}%
\mathbf{1}_{\left\{ X_{n}^{\left( N\right) }\left( t\right) =k,\text{ }%
J\left( t\right) =j\right\} },
\end{equation*}
where $\mathbf{1}_{\left\{ \bullet\right\} }$ is an indicative function.
Obviously, $Y_{k,j}^{\left( N\right) }\left( t\right) $ denotes the fraction
of stations with $k$ bikes and with the Markovian environment be at State $j$
at time $t$. It is easy to see that for $-L\leq k\leq K+L,$%
\begin{equation*}
0\leq Y_{k,j}^{\left( N\right) }\left( t\right)
\leq\sum\limits_{j=1}^{m}Y_{k,j}^{\left( N\right) }\left( t\right)
\leq\sum_{k=-L}^{K+L}\sum\limits_{j=1}^{m}Y_{k,j}^{\left( N\right) }\left(
t\right) =1.
\end{equation*}
Let%
\begin{equation*}
\mathbf{Y}_{k}^{\left( N\right) }\left( t\right) =\left( Y_{k,1}^{\left(
N\right) }\left( t\right) ,Y_{k,2}^{\left( N\right) }\left( t\right)
,\ldots,Y_{k,m}^{\left( N\right) }\left( t\right) \right)
\end{equation*}
and%
\begin{equation*}
\mathbf{Y}^{\left( N\right) }\left( t\right) =\left( \mathbf{Y}_{-L}^{\left(
N\right) }\left( t\right) ,\mathbf{Y}_{-L+1}^{\left( N\right) }\left(
t\right) ,\ldots,\mathbf{Y}_{K+L-1}^{\left( N\right) }\left( t\right) ,%
\mathbf{Y}_{K+L}^{\left( N\right) }\left( t\right) \right) ,
\end{equation*}
which is a row vector of size $\left( K+2L+1\right) m$. Then $\left\{
\mathbf{Y}^{\left( N\right) }\left( t\right) :t\geq0\right\} $ is an
empirical measure Markov process whose state space is given by $\Omega=\left[
0,1\right] ^{\left( K+2L+1\right) m}$.

To consider the empirical measure Markov process $\left\{ \mathbf{Y}^{\left(
N\right) }\left( t\right) :t\geq 0\right\} $, we write%
\begin{equation}
y_{k,j}^{\left( N\right) }\left( t\right) =E\left[ Y_{k,j}^{\left( N\right)
}\left( t\right) \right] ,\text{ \ \ }-L\leq k\leq K+L,1\leq j\leq m,
\label{Prob}
\end{equation}%
and%
\begin{equation*}
\mathbf{y}_{k}^{\left( N\right) }\left( t\right) =\left( y_{k,1}^{\left(
N\right) }\left( t\right) ,y_{k,2}^{\left( N\right) }\left( t\right) ,\ldots
,y_{k,m}^{\left( N\right) }\left( t\right) \right) ,
\end{equation*}%
\begin{equation*}
\mathbf{y}^{\left( N\right) }\left( t\right) =\left( \mathbf{y}_{-L}^{\left(
N\right) }\left( t\right) ,\mathbf{y}_{-L+1}^{\left( N\right) }\left(
t\right) ,\ldots ,\mathbf{y}_{K+L-1}^{\left( N\right) }\left( t\right) ,%
\mathbf{y}_{K+L}^{\left( N\right) }\left( t\right) \right) .
\end{equation*}%
In what follows we will apply the mean-field theory to set up a
block-structured system of mean-field equations whose purpose is to be able
to numerically compute the key vector $\mathbf{y}^{\left( N\right) }\left(
t\right) $.

\subsection{ A time-inhomogeneous MAP$\left( t\right) $/MAP$\left( t\right) $%
/1/K+2L+1 queue}

In the bike sharing system with $N$ identical stations and with Markovian
environment, we define $\mathbf{Q}^{\left( N\right) }\left( t\right) $ as
the number of bikes in a tagged station at time $t$. It is easy to see that
if an outside customer arrives at the tagged station and rents a bike, then $%
\mathbf{Q}^{\left( N\right) }\left( t\right) $ decreases by one; while if a
customer finishes his trip and returns a bike at the tagged station, $%
\mathbf{Q}^{\left( N\right) }\left( t\right) $ increases by one. Based on
this, we can understand that the Markov process $\left\{ \left( \mathbf{Q}%
^{\left( N\right) }\left( t\right) ,J\left( t\right) \right) :t\geq
0\right\} $ is a QBD process, which is further shown to well correspond to a
time-inhomogeneous MAP$\left( t\right) $/MAP$\left( t\right) $/1/K+2L+1
queue, where MAP$\left( t\right) $ is an instantaneous Markov arrival
process with a matrix descriptor $\left( C\left( t\right) ,D\left( t\right)
\right) $ of size $m$, e.g., see Subsections 8.2.5 and 8.2.6 of Chapter 8 in
Li \cite{Li:2010} for more details.

In the MAP$\left( t\right) $/MAP$\left( t\right) $/1/K+2L+1 queue, it is
easy to understand that the customers are virtual from the bikes, that is,
the bikes are the virtual customers. Thus a virtual customer's arrival is a
bike returned to the tagged station; while a virtual customer's service
completion is a bike rented from the tagged station. Thus, here we call
arrival (or service) to be virtual arrival (or virtual service).

The following theorem provides expressions for the instantaneous virtual
arrival rate $\xi _{l,j}^{\left( N\right) }\left( t\right) $ and the
instantaneous virtual service rate $\eta _{k,j}^{\left( N\right) }\left(
t\right) $ in this time-inhomogeneous queueing system. Note that the two
instantaneous rates play a key role in our later study.

\begin{The}
\label{The:Queue}In the time-inhomogeneous MAP$\left( t\right) $/MAP$\left(
t\right) $/1/K+2L+1 queue, we have

(a) the instantaneous virtual service rate is given by%
\begin{equation}
\eta _{k,j}^{\left( N\right) }\left( t\right) =\left\{
\begin{array}{l}
\lambda _{j},\text{ \ \ \ }1\leq k\leq K+L,\text{ \ \ \ \ }1\leq j\leq m, \\
\lambda _{j}\alpha ,\text{ \ }-\left( L-1\right) \leq k\leq 0,\text{ \ }%
1\leq j\leq m,%
\end{array}%
\right.  \label{eq-1}
\end{equation}%
which is independent of the number $N$.

(b) For $1\leq j\leq m$, the instantaneous virtual arrival rate is given by%
\begin{equation}
\xi _{l,j}^{\left( N\right) }\left( t\right) =\left\{
\begin{array}{ll}
\frac{\mu _{j}}{N}\left\{ C+\left( N-1\right) \left[ C-\sum%
\limits_{k=1}^{K+L}ky_{k,j}^{\left( N\right) }\left( t\right) \right. \right.
&  \\
\left. \left. +\sum\limits_{k=K}^{K+L-1}\frac{\left( 1-\beta \right)
y_{k,j}^{\left( N\right) }\left( t\right) }{\left[ 1-\left( 1-\beta \right)
y_{k,j}^{\left( N\right) }\left( t\right) \right] ^{2}}+\frac{%
y_{K+L,j}^{\left( N\right) }\left( t\right) }{\left[ 1-y_{K+L,j}^{\left(
N\right) }\left( t\right) \right] ^{2}}\right] \right\} , & \text{\ }-L\leq
l\leq 0, \\
\frac{\mu _{j}}{N}\left\{ C-l+\left( N-1\right) \left[ C-\sum%
\limits_{k=1}^{K+L}ky_{k,j}^{\left( N\right) }\left( t\right) \right. \right.
&  \\
\left. \left. +\sum\limits_{k=K}^{K+L-1}\frac{\left( 1-\beta \right)
y_{k,j}^{\left( N\right) }\left( t\right) }{\left[ 1-\left( 1-\beta \right)
y_{k,j}^{\left( N\right) }\left( t\right) \right] ^{2}}+\frac{%
y_{K+L,j}^{\left( N\right) }\left( t\right) }{\left[ 1-y_{K+L,j}^{\left(
N\right) }\left( t\right) \right] ^{2}}\right] \right\} , & 1\leq l\leq C-1,
\\
\frac{\mu _{j}}{N}\left\{ \left( N-1\right) \left[ C-\sum%
\limits_{k=1}^{K+L}ky_{k,j}^{\left( N\right) }\left( t\right) \right. \right.
&  \\
\left. \left. +\sum\limits_{k=K}^{K+L-1}\frac{\left( 1-\beta \right)
y_{k,j}^{\left( N\right) }\left( t\right) }{\left[ 1-\left( 1-\beta \right)
y_{k,j}^{\left( N\right) }\left( t\right) \right] ^{2}}+\frac{%
y_{K+L,j}^{\left( N\right) }\left( t\right) }{\left[ 1-y_{K+L,j}^{\left(
N\right) }\left( t\right) \right] ^{2}}\right] \right\} , & C\leq l\leq K-1,
\\
\beta \frac{\mu _{j}}{N}\left\{ \left( N-1\right) \left[ C-\sum%
\limits_{k=1}^{K+L}ky_{k,j}^{\left( N\right) }\left( t\right) \right. \right.
&  \\
\left. \left. +\sum\limits_{k=K}^{K+L-1}\frac{\left( 1-\beta \right)
y_{k,j}^{\left( N\right) }\left( t\right) }{\left[ 1-\left( 1-\beta \right)
y_{k,j}^{\left( N\right) }\left( t\right) \right] ^{2}}+\frac{%
y_{K+L,j}^{\left( N\right) }\left( t\right) }{\left[ 1-y_{K+L,j}^{\left(
N\right) }\left( t\right) \right] ^{2}}\right] \right\} , & K\leq l\leq
K+L-1.%
\end{array}%
\right.  \label{eq-2}
\end{equation}
\end{The}

\textbf{Proof:} \underline{The proof of (2)}. When a customer arrives at a
tagged station, there exist two different cases:

\textbf{Case (a)} If the station has at least one bike (that is, $1\leq
k\leq K+L$), then he immediately rents a bike and leaves the station, that
is, the virtual service is completed. Hence $\eta_{k,j}^{\left( N\right)
}\left( t\right) =$ $\lambda_{j}$ for $1\leq j\leq m$.

\textbf{Case (b)} If the station has no bike (that is, $-L+1\leq k\leq 0$),
then he has two choices: he directly leaves this system with the probability
$1-\alpha $; or he enters a waiting place with the probability $\alpha $ in
order to wait for renting a future available bike. Clearly, the virtual
service has not been completed yet but this also leads to the shortage of
virtual customers (or bikes). In this case, the rate $\lambda _{j}\alpha $
expresses the transition speed of that the number of bikes at the tagged
station from $k$ to $k-1$ for $-\left( L-1\right) \leq k\leq 0$. Thus $\eta
_{k,j}^{\left( N\right) }\left( t\right) =\lambda _{j}\alpha $ for $1\leq
j\leq m$, and it is independent of the number $k=-\left( L-1\right) ,-\left(
L-2\right) ,\ldots ,1,0$.

Based on Cases (a) and (b), when the Markovian environment $J\left( t\right)
=j$, we have%
\begin{equation*}
\eta _{k,j}^{\left( N\right) }\left( t\right) =\left\{
\begin{array}{l}
\lambda _{j},\text{ \ \ \ \ }1\leq k\leq K+L,\text{ \ \ \ \ }1\leq j\leq m,
\\
\lambda _{j}\alpha ,\text{ \ \ }-\left( L-1\right) \leq k\leq 0,\text{\ \ }%
1\leq j\leq m,%
\end{array}%
\right.
\end{equation*}%
which is independent of $k=-\left( L-1\right) ,-\left( L-2\right) ,\ldots
,K+L-1,K+L$.

\underline{The proof of (3)}. The proof of (3) is a bit complicated due to
applications of the mean-field theory. Note that Figure 3 describes the
state transitions of the process that the bikes are returned at the tagged
station. Since the bikes can not be returned to a full station, either the
user enters a waiting place in order to wait for a future available parking
place with probability $\beta $, or he immediately rides his bike to another
station to find an available parking place with probability $1-\beta $.
Based on this, we give the instantaneous virtual arrival rate $\xi
_{l,j}^{\left( N\right) }\left( t\right) $ by means of a
probability-analytic method as follows:%
\begin{equation*}
\xi _{l,j}^{\left( N\right) }\left( t\right) =\frac{1}{N}\cdot \mu _{j}\cdot
\text{the number of bikes ridden on all the roads at time $t$}.
\end{equation*}%
Note that the number of bikes ridden on all the roads contains two parts:
(i) The number $n_{1}$ of bikes ridden from the tagged station is given by%
\begin{equation*}
n_{1}=\left\{
\begin{array}{ll}
C, & -L\leq l\leq 0, \\
C-l, & 1\leq l\leq C-1, \\
0, & C\leq l\leq K+L-1,%
\end{array}%
\right.
\end{equation*}%
and (ii) the number of bikes having been ridden from the other $N-1$ station
is given by%
\begin{align*}
& \left( N-1\right) \left[ \text{\textit{the average number of bikes
directly ridden from one of the $N-1$ stations}}\right. \\
& \left. +\text{\textit{the average number of bikes which can not be
returned to a full station with at least}}\right. \\
& \left. \text{\textit{\ two retries}}\right] .
\end{align*}%
By using the mean-field theory, the average number of bikes ridden from the
tagged station is given by $C-\sum_{k=1}^{K+L}ky_{k,j}^{\left( N\right)
}\left( t\right) $. While the average number of bikes which can not be
returned to a full station with at least two retries is given a detailed
computation below.

\begin{figure}[ptb]
\setlength{\abovecaptionskip}{0cm} \setlength{\belowcaptionskip}{-0cm} %
\centering              \includegraphics[width=10cm]{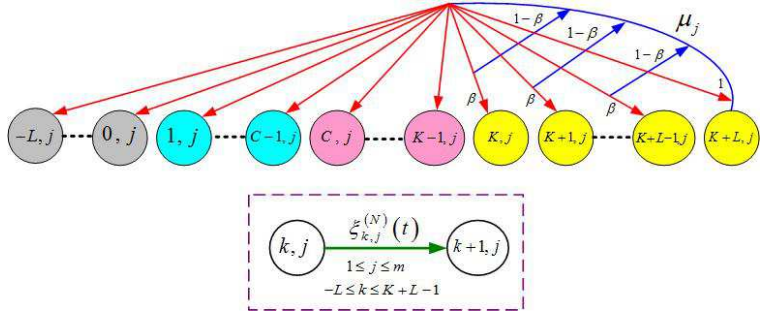} \newline
\caption{The state transition relation of the queueing process in the tagged
station}
\label{figure:figure-2}
\end{figure}

Based on the above analysis, our computation for deriving the instantaneous
virtual arrival rate $\xi_{l,j}^{\left( N\right) }\left( t\right) $ is
divided into the following four cases.

\textbf{Case (a):} When $-L\leq l\leq 0$, $1\leq j\leq m$, we need to study
three different classes for the initial distribution of bikes in the tagged
station. Note that in the last two classes, users who arrive at a full
station can not return their bikes at the station.

Class-1: The initial $C$ bikes in the tagged station are all rented on the
roads. Using the mean-field theory, we get that the average number of bikes
rented on the road from the other $N-1$ stations is given by%
\begin{equation*}
\left( N-1\right) \left[ C-\sum_{k=1}^{K+L}ky_{k,j}^{\left( N\right) }\left(
t\right) \right] ,
\end{equation*}%
where $\sum\nolimits_{k=1}^{K+L}ky_{k,j}^{\left( N\right) }\left( t\right) $
is the average number of bikes parked in the tagged station. Thus, the
average number of bikes rented on the roads from the $N$ stations is given by%
\begin{equation*}
C+\left( N-1\right) \left[ C-\sum_{k=1}^{K+L}ky_{k,j}^{\left( N\right)
}\left( t\right) \right] .
\end{equation*}

Class-2: A customer\ finishes his trip and arrives at a tagged station in
which there are $k$ bikes for $K\leq k\leq K+L-1.$ It is clear that the
tagged station is full so that the customer has to re-ride the bike in order
to return the bike to another station with the probability $1-\beta $. The
average number of such re-riding bike is given by%
\begin{align*}
& \sum\limits_{k=K}^{K+L-1}\left\{ \left( 1-\beta \right) y_{k,j}^{\left(
N\right) }\left( t\right) +2\left[ \left( 1-\beta \right) y_{k,j}^{\left(
N\right) }\left( t\right) \right] ^{2}+3\left[ \left( 1-\beta \right)
y_{k,j}^{\left( N\right) }\left( t\right) \right] ^{3}+\cdots \right\} \\
& =\sum\limits_{k=K}^{K+L-1}\frac{\left( 1-\beta \right) y_{k,j}^{\left(
N\right) }\left( t\right) }{\left[ 1-\left( 1-\beta \right) y_{k,j}^{\left(
N\right) }\left( t\right) \right] ^{2}},
\end{align*}%
where $n\left[ \left( 1-\beta \right) y_{k,j}^{\left( N\right) }\left(
t\right) \right] ^{n}$ is the average number of re-riding bikes of $n$
customers, and $x+2x^{2}+3x^{3}+\cdots =x/\left( 1-x\right) ^{2}$.

Class-3: A customer finishes his trip and arrives at a tagged station in
which there are $K+L$ bikes. In this case, there is neither a bike-parking
place nor a user-waiting place, hence the customer has to re-ride the bike
in order to return the bike at another station with the probability $1$. The
average number of such re-riding bikes is given by%
\begin{equation*}
y_{K+L,j}^{\left( N\right) }\left( t\right) +2\left[ y_{K+L,j}^{\left(
N\right) }\left( t\right) \right] ^{2}+3\left[ y_{K+L,j}^{\left( N\right)
}\left( t\right) \right] ^{3}+\cdots =\frac{y_{K+L,j}^{\left( N\right)
}\left( t\right) }{\left[ 1-y_{K+L,j}^{\left( N\right) }\left( t\right) %
\right] ^{2}}.
\end{equation*}

Summarizing the above analysis, the instantaneous virtual arrival rate is
given by{\small
\begin{equation*}
\xi_{l,j}^{\left( N\right) }\left( t\right) =\frac{\mu_{j}}{N}\left\{
C+\left( N-1\right) \left[ C-\sum\limits_{k=1}^{K+L}ky_{k,j}^{\left(
N\right) }\left( t\right) +\sum\limits_{k=K}^{K+L-1}\frac{\left(
1-\beta\right) y_{k,j}^{\left( N\right) }\left( t\right) }{\left[ 1-\left(
1-\beta\right) y_{k,j}^{\left( N\right) }\left( t\right) \right] ^{2}}+\frac{%
y_{K+L,j}^{\left( N\right) }\left( t\right) }{\left[ 1-y_{K+L,j}^{\left(
N\right) }\left( t\right) \right] ^{2}}\right] \right\} .
\end{equation*}
}

\textbf{Case (b):} When $1\leq l\leq C-1$, the only difference of our
derivation from Case (a) is to replace the initial $C$ bikes by the initial $%
C-l$ bikes in the tagged station. Thus we get{\small
\begin{equation*}
\xi_{l,j}^{\left( N\right) }\left( t\right) =\frac{\mu_{j}}{N}\left\{
C-l+\left( N-1\right) \left[ C-\sum\limits_{k=1}^{K+L}ky_{k,j}^{\left(
N\right) }\left( t\right) +\sum\limits_{k=K}^{K+L-1}\frac{\left(
1-\beta\right) y_{k,j}^{\left( N\right) }\left( t\right) }{\left[ 1-\left(
1-\beta\right) y_{k,j}^{\left( N\right) }\left( t\right) \right] ^{2}}+\frac{%
y_{K+L,j}^{\left( N\right) }\left( t\right) }{\left[ 1-y_{K+L,j}^{\left(
N\right) }\left( t\right) \right] ^{2}}\right] \right\} .
\end{equation*}
}

\textbf{Case (c):} When $C\leq l\leq K-1$, the only difference of our
derivation from Case (a) is that the initial $C$ bikes in this station are
all parked in the tagged station. Hence we obtain%
\begin{equation*}
\xi_{l,j}^{\left( N\right) }\left( t\right) =\frac{\mu_{j}}{N}\left\{ \left(
N-1\right) \left[ C-\sum\limits_{k=1}^{K+L}ky_{k,j}^{\left( N\right) }\left(
t\right) +\sum\limits_{k=K}^{K+L-1}\frac{\left( 1-\beta\right)
y_{k,j}^{\left( N\right) }\left( t\right) }{\left[ 1-\left( 1-\beta\right)
y_{k,j}^{\left( N\right) }\left( t\right) \right] ^{2}}+\frac{%
y_{K+L,j}^{\left( N\right) }\left( t\right) }{\left[ 1-y_{K+L,j}^{\left(
N\right) }\left( t\right) \right] ^{2}}\right] \right\} .
\end{equation*}

\textbf{Case (d):} When $K\leq l\leq K+L-1$, the only difference of our
derivation from Case (c) is that when a customer finishes his trip and
arrives at the tagged station, he enters the waiting places in order to wait
for an empty parking place with the probability $\beta$. This gives%
\begin{equation*}
\xi_{l,j}^{\left( N\right) }\left( t\right) =\beta\frac{\mu_{j}}{N}\left\{
\left( N-1\right) \left[ C-\sum\limits_{k=1}^{K+L}ky_{k,j}^{\left( N\right)
}\left( t\right) +\sum\limits_{k=K}^{K+L-1}\frac{\left( 1-\beta\right)
y_{k,j}^{\left( N\right) }\left( t\right) }{\left[ 1-\left( 1-\beta\right)
y_{k,j}^{\left( N\right) }\left( t\right) \right] ^{2}}+\frac{%
y_{K+L,j}^{\left( N\right) }\left( t\right) }{\left[ 1-y_{K+L,j}^{\left(
N\right) }\left( t\right) \right] ^{2}}\right] \right\} .
\end{equation*}

Summarizing the above four cases, we obtain all the expressions given in
(3). This completes the proof. \textbf{{\rule{0.08in}{0.08in}}}

\subsection{A block-structured system of mean-field equations}

Based on the time-inhomogeneous MAP$\left( t\right) $/MAP$\left( t\right) $%
/1/K+2L+1 queue, it is convenient to describe the time-inhomogeneous QBD
process $\left\{ \left( \mathbf{Q}^{\left( N\right) }\left( t\right)
,J\left( t\right) \right) :t\geq 0\right\} $ whose state transition relation
is depicted in Figure 4. At the same time, a useful relation related to (1)
for understanding the state probability distribution of the QBD process is
given by

\begin{figure}[ptb]
\setlength{\abovecaptionskip}{0.cm} \setlength{\belowcaptionskip}{-0.cm} %
\centering              \includegraphics[width=10cm]{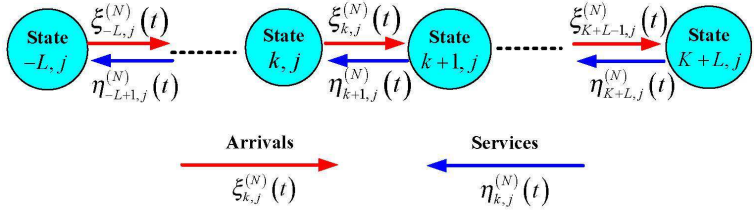} \newline
\caption{The state transitions of the time-inhomogeneous QBD process}
\label{figure:figure-3}
\end{figure}

\begin{equation*}
y_{k,j}^{\left( N\right) }\left( t\right) =P\left\{ \mathbf{Q}^{\left(
N\right) }\left( t\right) =k,J\left( t\right) =j\right\} ,\text{\ }-L\leq
k\leq K+L,1\leq j\leq m.
\end{equation*}%
Based on the instantaneous virtual arrival rate $\xi _{l,j}^{\left( N\right)
}\left( t\right) $, we set that for $-L\leq k\leq K+L-1$
\begin{equation*}
\Psi _{k}^{\left( N\right) }\left( t\right) =\left(
\begin{array}{cccc}
0 & \xi _{k,1}^{\left( N\right) }\left( t\right) w_{1,2} & \cdots & \xi
_{k,1}^{\left( N\right) }\left( t\right) w_{1,m} \\
\xi _{k,2}^{\left( N\right) }\left( t\right) w_{2,1} & 0 & \cdots & \xi
_{k,2}^{\left( N\right) }\left( t\right) w_{2,m} \\
\vdots & \vdots & \ddots & \vdots \\
\xi _{k,m}^{\left( N\right) }\left( t\right) w_{m,1} & \xi _{k,m}^{\left(
N\right) }\left( t\right) w_{m,2} & \cdots & 0%
\end{array}%
\right) ,
\end{equation*}%
\begin{equation*}
\widehat{\Psi }_{k}^{\left( N\right) }\left( t\right) =\text{diag}\left( \xi
_{k,1}^{\left( N\right) }\left( t\right) w_{1,1},\xi _{k,2}^{\left( N\right)
}\left( t\right) w_{2,2},\ldots ,\xi _{k,m}^{\left( N\right) }\left(
t\right) w_{m,m}\right) .
\end{equation*}%
Similarly, it is easy to see from the instantaneous virtual service rate $%
\eta _{k,j}^{\left( N\right) }\left( t\right) $\ that for $-L+1\leq k\leq 0$%
\begin{equation*}
\Phi _{k}^{\left( N\right) }\left( t\right) =\left(
\begin{array}{cccc}
0 & \lambda _{1}\alpha w_{1,2} & \cdots & \lambda _{1}\alpha w_{1,m} \\
\lambda _{2}\alpha w_{2,1} & 0 & \cdots & \lambda _{2}\alpha w_{2,m} \\
\vdots & \vdots & \ddots & \vdots \\
\lambda _{m}\alpha w_{m,1} & \lambda _{m}\alpha w_{m,2} & \cdots & 0%
\end{array}%
\right) \overset{\text{Def}}{=}\Phi \left( \alpha \right) ,
\end{equation*}%
\begin{equation*}
\widehat{\Phi }_{k}^{\left( N\right) }\left( t\right) =\text{diag}\left(
\lambda _{1}\alpha w_{1,1},\lambda _{2}\alpha w_{2,2},\ldots ,\lambda
_{m}\alpha w_{m,m}\right) \overset{\text{Def}}{=}\widehat{\Phi }\left(
\alpha \right) ;
\end{equation*}%
while for $1\leq k\leq K+L$%
\begin{equation*}
\Phi _{k}^{\left( N\right) }\left( t\right) =\Phi \left( 1\right)
\end{equation*}%
and%
\begin{equation*}
\widehat{\Phi }_{k}^{\left( N\right) }\left( t\right) =\widehat{\Phi }\left(
1\right) ,
\end{equation*}%
both of which are due to $\alpha =1$ for $1\leq k\leq K+L$.

For the time-inhomogeneous QBD process $\left\{ \left( \mathbf{Q}^{\left(
N\right) }\left( t\right) ,J\left( t\right) \right) :t\geq 0\right\} $, it
follows from Figure 4 that the vector $\mathbf{y}^{\left( N\right) }\left(
t\right) =\left( \mathbf{y}_{-L}^{\left( N\right) }\left( t\right) ,\mathbf{y%
}_{-L+1}^{\left( N\right) }\left( t\right) ,\ldots ,\mathbf{y}%
_{K+L-1}^{\left( N\right) }\left( t\right) ,\mathbf{y}_{K+L}^{\left(
N\right) }\left( t\right) \right) $ satisfies a block-structured system of
mean-field (or ordinary differential) equations as follows:
\begin{equation}
\frac{\text{d}}{\text{d}t}\mathbf{y}_{-L}^{\left( N\right) }\left( t\right) =%
\mathbf{y}_{-L}^{\left( N\right) }\left( t\right) \widehat{\Psi }%
_{-L}^{\left( N\right) }\left( t\right) +\mathbf{y}_{-L+1}^{\left( N\right)
}\left( t\right) \Phi _{-L+1}^{\left( N\right) }\left( t\right) ,
\label{eq-3}
\end{equation}%
for $-L+1\leq k\leq K+L-1$%
\begin{equation}
\frac{\text{d}}{\text{d}t}\mathbf{y}_{k}^{\left( N\right) }\left( t\right) =%
\mathbf{y}_{k-1}^{\left( N\right) }\left( t\right) \Psi _{k-1}^{\left(
N\right) }\left( t\right) +\mathbf{y}_{k}^{\left( N\right) }\left( t\right) %
\left[ \widehat{\Psi }_{k}^{\left( N\right) }\left( t\right) +\widehat{\Phi }%
_{k}^{\left( N\right) }\left( t\right) \right] +\mathbf{y}_{k+1}^{\left(
N\right) }\left( t\right) \Phi _{k+1}^{\left( N\right) }\left( t\right) ,
\label{eq-4}
\end{equation}%
\begin{equation}
\frac{\text{d}}{\text{d}t}\mathbf{y}_{K+L}^{\left( N\right) }\left( t\right)
=\mathbf{y}_{K+L}^{\left( N\right) }\left( t\right) \widehat{\Phi }%
_{K+L}^{\left( N\right) }\left( t\right) +\mathbf{y}_{K+L-1}^{\left(
N\right) }\left( t\right) \Psi _{K+L-1}^{\left( N\right) }\left( t\right) ,
\label{eq-5}
\end{equation}%
with the boundary condition%
\begin{equation}
\sum_{k=-L}^{K+L}\mathbf{y}_{k}^{\left( N\right) }\left( t\right) e=1,
\label{eq-6-1}
\end{equation}%
and with the initial condition
\begin{equation}
\mathbf{y}_{k}^{\left( N\right) }\left( 0\right) =g_{k},\text{ \ }-L\leq
k\leq K+L,  \label{eq-6-2}
\end{equation}%
and%
\begin{equation*}
g_{k}=\left( g_{k,1},g_{k,2},\ldots ,g_{k,m}\right) ,
\end{equation*}%
\begin{equation*}
\mathbf{g}=\left( g_{-L},g_{-\left( L-1\right) },\cdots
,g_{K+L-1},g_{K+L}\right)
\end{equation*}%
is a probability vector of size $\left( K+2L+1\right) m$.

For convenience of description, we write the mean-field equations (\ref{eq-3}%
) to (\ref{eq-6-2}) into a matrix version as follows:%
\begin{equation}
\frac{\text{d}}{\text{d}t}\mathbf{y}^{\left( N\right) }\left( t\right) =%
\mathbf{y}^{\left( N\right) }\left( t\right) \mathbf{V}_{\mathbf{y}^{\left(
N\right) }\left( t\right) },  \label{eq-7}
\end{equation}
with the boundary and initial conditions%
\begin{equation}
\mathbf{y}^{\left( N\right) }\left( t\right) e=1,\text{ \ }\mathbf{y}%
^{\left( N\right) }\left( 0\right) =\mathbf{g,}  \label{eq-7-1}
\end{equation}
where
\begin{equation}
\mathbf{V}_{\mathbf{y}^{\left( N\right) }\left( t\right) }=\left(
\begin{array}{cccc}
A_{1,1} & A_{1,2} &  &  \\
A_{2,1} & A_{2,2} & A_{2,3} &  \\
& A_{3,2} & A_{3,3} & A_{3,4} \\
&  & A_{4,3} & A_{4,4}%
\end{array}
\right) ,  \label{eq-7-2}
\end{equation}%
\begin{equation*}
\Delta_{k}^{\left( N\right) }\left( t\right) =\widehat{\Phi}_{k}^{\left(
N\right) }\left( t\right) +\widehat{\Psi}_{k}^{\left( N\right) }\left(
t\right) ,\text{ \ }-\left( L-1\right) \leq k\leq K+L-1,
\end{equation*}%
\begin{equation*}
A_{1,1}=\left(
\begin{array}{ccccc}
\widehat{\Psi}_{-L}^{\left( N\right) }\left( t\right) & \Psi_{-L}^{\left(
N\right) }\left( t\right) &  &  &  \\
\Phi_{-\left( L-1\right) }^{\left( N\right) }\left( t\right) &
\Delta_{-\left( L-1\right) }^{\left( N\right) }\left( t\right) &
\Psi_{-\left( L-1\right) }^{\left( N\right) }\left( t\right) &  &  \\
& \ddots & \ddots & \ddots &  \\
&  & \Phi_{-1}^{\left( N\right) }\left( t\right) & \Delta_{-1}^{\left(
N\right) }\left( t\right) & \Psi_{-1}^{\left( N\right) }\left( t\right) \\
&  &  & \Phi_{0}^{\left( N\right) }\left( t\right) & \Delta_{0}^{\left(
N\right) }\left( t\right)%
\end{array}
\right) ,
\end{equation*}%
\begin{equation*}
A_{1,2}=\left(
\begin{array}{cccc}
&  &  &  \\
&  &  &  \\
&  &  &  \\
\Psi_{0}^{\left( N\right) }\left( t\right) &  &  &
\end{array}
\right) ,\text{ \ }A_{2,1}=\left(
\begin{array}{cccc}
&  &  & \Phi_{1}^{\left( N\right) }\left( t\right) \\
&  &  &  \\
&  &  &  \\
&  &  &
\end{array}
\right) ,
\end{equation*}%
\begin{equation*}
A_{2,2}=\left(
\begin{array}{ccccc}
\Delta_{1}^{\left( N\right) }\left( t\right) & \Psi_{1}^{\left( N\right)
}\left( t\right) &  &  &  \\
\Phi_{2}^{\left( N\right) }\left( t\right) & \Delta_{2}^{\left( N\right)
}\left( t\right) & \Psi_{2}^{\left( N\right) }\left( t\right) &  &  \\
& \ddots & \ddots & \ddots &  \\
&  & \Phi_{L-2}^{\left( N\right) }\left( t\right) & \Delta_{L-2}^{\left(
N\right) }\left( t\right) & \Psi_{L-2}^{\left( N\right) }\left( t\right) \\
&  &  & \Phi_{L-1}^{\left( N\right) }\left( t\right) & \Delta _{L-1}^{\left(
N\right) }\left( t\right)%
\end{array}
\right) ,
\end{equation*}%
\begin{equation*}
A_{2,3}=\left(
\begin{array}{cccc}
&  &  &  \\
&  &  &  \\
&  &  &  \\
\Psi_{L-1}^{\left( N\right) }\left( t\right) &  &  &
\end{array}
\right) ,\text{ \ }A_{3,2}=\left(
\begin{array}{cccc}
&  &  & \Phi_{C}^{\left( N\right) }\left( t\right) \\
&  &  &  \\
&  &  &  \\
&  &  &
\end{array}
\right) ,
\end{equation*}%
\begin{equation*}
A_{3,3}=\left(
\begin{array}{ccccc}
\Delta_{C}^{\left( N\right) }\left( t\right) & \Psi_{C}^{\left( N\right)
}\left( t\right) &  &  &  \\
\Phi_{C+1}^{\left( N\right) }\left( t\right) & \Delta_{C+1}^{\left( N\right)
}\left( t\right) & \Psi_{C+1}^{\left( N\right) }\left( t\right) &  &  \\
& \ddots & \ddots & \ddots &  \\
&  & \Phi_{K-2}^{\left( N\right) }\left( t\right) & \Delta_{K-2}^{\left(
N\right) }\left( t\right) & \Psi_{K-2}^{\left( N\right) }\left( t\right) \\
&  &  & \Phi_{K-1}^{\left( N\right) }\left( t\right) & \Delta _{K-1}^{\left(
N\right) }\left( t\right)%
\end{array}
\right) ,
\end{equation*}%
\begin{equation*}
A_{3,4}=\left(
\begin{array}{cccc}
&  &  &  \\
&  &  &  \\
&  &  &  \\
\Psi_{K-1}^{\left( N\right) }\left( t\right) &  &  &
\end{array}
\right) ,\text{ \ }A_{4,3}=\left(
\begin{array}{cccc}
&  &  & \Phi_{K}^{\left( N\right) }\left( t\right) \\
&  &  &  \\
&  &  &  \\
&  &  &
\end{array}
\right) ,
\end{equation*}%
\begin{equation*}
A_{4,4}=\left(
\begin{array}{ccccc}
\Delta_{K}^{\left( N\right) }\left( t\right) & \Psi_{K}^{\left( N\right)
}\left( t\right) &  &  &  \\
\Phi_{K+1}^{\left( N\right) }\left( t\right) & \Delta_{K+1}^{\left( N\right)
}\left( t\right) & \Psi_{K+1}^{\left( N\right) }\left( t\right) &  &  \\
& \ddots & \ddots & \ddots &  \\
&  & \Phi_{K+L-1}^{\left( N\right) }\left( t\right) & \Delta
_{K+L-1}^{\left( N\right) }\left( t\right) & \Psi_{K+L-1}^{\left( N\right)
}\left( t\right) \\
&  &  & \Phi_{K+L}^{\left( N\right) }\left( t\right) & \widehat{\Phi }%
_{K+L}^{\left( N\right) }\left( t\right)%
\end{array}
\right) .
\end{equation*}

\begin{Rem}
To set up the block-structured system of mean-field equations, it is a key
to observe two figures: Figure 1 shows all the original parameters of the
queueing process under Markovian environment when one isolated station is
paid attention to. Figure 4 further gives mean-field expressions for the
transition rates of the QBD process between two neighboring levels by
considering the weak interactions among the $N$ stations of the bike sharing
system in terms of the mean-field theory.
\end{Rem}

\section{The Martingale Limits}

In this section, we apply the martingale limit theory to prove the
asymptotic independence of this bike sharing system, that is, the sequence $%
\left\{ \mathbf{Y}^{\left( N\right) }\left( t\right) ,t\geq0\right\} $ of
Markov processes asymptotically approaches a single trajectory identified by
a solution to the block-structured system of limiting mean-field equations.

For the vector $\mathbf{g}=\left( g_{-L},g_{-\left( L-1\right) },\ldots
,g_{K+L-1},g_{K+L}\right) $ where $g_{k}=\left( g_{k,1},g_{k,2},\ldots
,g_{k,m}\right) $, we set%
\begin{equation*}
\Omega _{N}=\left\{ \mathbf{g}:\mathbf{g}\geq 0,\text{ }\mathbf{g}e=1,\text{
}N\mathbf{g}\text{ is a vector of nonnegative integers}\right\}
\end{equation*}%
and%
\begin{equation*}
\Omega =\left\{ \mathbf{g}:\mathbf{g}\geq 0,\text{ }\mathbf{g}e=1\right\} .
\end{equation*}%
Obviously, $\Omega _{N}\subset \Omega $. In the vector space $\Omega $ (or $%
\Omega _{N}$), we take a metric%
\begin{equation*}
\rho \left( \mathbf{g},\mathbf{g}^{^{\prime }}\right) =\max_{-L\leq k\leq
K+L}\max_{1\leq j\leq m}\left\{ |g_{k,j}-g_{k,j}^{\prime }|\right\} ,\text{
\ \ }\mathbf{g},\mathbf{g}^{\prime }\in \Omega .
\end{equation*}%
Note that under the metric $\rho \left( \mathbf{g},\mathbf{g}^{^{\prime
}}\right) $, the vector space $\Omega $ (or $\Omega _{N}$) is separable and
compact.

Now we consider the Markov process $\left\{ \mathbf{Y}^{\left( N\right)
}\left( t\right) ,t\geq 0\right\} $ on state space $\Omega _{N}$ for $%
N=1,2,3,\ldots $. Note that the stochastic evolution of this bike sharing
system is described as the Markov process $\left\{ \mathbf{Y}^{\left(
N\right) }\left( t\right) ,t\geq 0\right\} $, and%
\begin{equation*}
\frac{\text{d}}{\text{d}t}\mathbf{Y}^{(N)}(t)=\mathbf{A}_{N}\text{ }f(%
\mathbf{Y}^{(N)}(t)),
\end{equation*}%
where $\mathbf{A}_{N}$ acting on functions $f:\Omega _{N}\rightarrow \mathbf{%
C}^{1}$ is the generating operator of the Markov process $\left\{ \mathbf{Y}%
^{\left( N\right) }\left( t\right) ,t\geq 0\right\} $, and%
\begin{equation}
\mathbf{A}_{N}=\mathbf{A}_{N}^{\text{renting}}+\mathbf{A}_{N}^{\text{%
returning}}+\mathbf{A}_{N}^{\text{environment}},  \label{oper-1}
\end{equation}%
where%
\begin{align*}
\mathbf{A}_{N}^{\text{renting}}f(\mathbf{g})=& N\sum_{j=1}^{m}\lambda
_{j}\sum_{k=1}^{K+L}g_{k,j}\left[ f\left( \mathbf{g}-\frac{e_{k,j}}{N}%
\right) -f\left( \mathbf{g}\right) \right] \\
& +N\alpha \sum_{j=1}^{m}\lambda _{j}\sum_{k=-\left( L-1\right) }^{0}g_{k,j}
\left[ f\left( \mathbf{g}-\frac{e_{k,j}}{N}\right) -f\left( \mathbf{g}%
\right) \right] ,
\end{align*}%
\begin{equation*}
\mathbf{A}_{N}^{\text{environment}}f(\mathbf{g})=N\sum_{i=1}^{m}%
\sum_{j=1}^{m}\sum_{k=-L}^{K+L}g_{k,i}w_{i,j}\left[ f\left( \mathbf{g}-\frac{%
e_{k,i}}{N}+\frac{e_{k,j}}{N}\right) -f\left( \mathbf{g}\right) \right] ,
\end{equation*}%
\begin{align*}
\mathbf{A}_{N}^{\text{returning}}f(\mathbf{g})=& \left\{
\sum_{j=1}^{m}\sum_{l=-L}^{0}\mu _{j}g_{l,j}\Theta ^{\left( N\right) }\left(
0\right) +\sum_{j=1}^{m}\sum_{l=1}^{C-1}\mu _{j}g_{l,j}\Theta ^{\left(
N\right) }\left( l\right) \right. \\
& \left. +\sum_{j=1}^{m}\sum_{l=C}^{K-1}\mu _{j}g_{l,j}\Theta ^{\left(
N\right) }\left( C\right) +\beta \sum_{j=1}^{m}\sum_{l=K}^{K+L-1}\mu
_{j}g_{l,j}\Theta ^{\left( N\right) }\left( C\right) \right\} \left[ f\left(
\mathbf{g}+\frac{e_{k,j}}{N}\right) -f\left( \mathbf{g}\right) \right] ,
\end{align*}%
and for $0\leq l\leq C$%
\begin{equation*}
\Theta ^{\left( N\right) }\left( l\right) =\left\{ C-l+\left( N-1\right)
\left[ C-\sum\limits_{k=1}^{K+L}kg_{k,j}+\sum\limits_{k=K}^{K+L-1}\frac{%
\left( 1-\beta \right) g_{k,j}}{\left[ 1-\left( 1-\beta \right) g_{k,j}%
\right] ^{2}}+\frac{g_{K+L,j}}{\left[ 1-g_{K+L,j}\right] ^{2}}\right]
\right\} .
\end{equation*}%
When $N\rightarrow \infty ,$ it is easy to check that%
\begin{equation*}
N\left[ f\left( \mathbf{g}+\frac{e_{k,j}}{N}\right) -f\left( \mathbf{g}%
\right) \right] \rightarrow \frac{\partial }{\partial g_{k,j}}f(\mathbf{g}),
\end{equation*}%
\begin{equation*}
N\left[ f\left( \mathbf{g}-\frac{e_{k,j}}{N}\right) -f\left( \mathbf{g}%
\right) \right] \rightarrow -\frac{\partial }{\partial g_{k,j}}f(\mathbf{g}),
\end{equation*}%
\begin{equation*}
\left[ f\left( \mathbf{g}-\frac{e_{k,i}}{N}+\frac{e_{k,j}}{N}\right)
-f\left( \mathbf{g}\right) \right] \rightarrow -\frac{\partial }{\partial
g_{k,i}}f(\mathbf{g})+\frac{\partial }{\partial g_{k,j}}f(\mathbf{g}),
\end{equation*}%
and for $0\leq l\leq C$%
\begin{align*}
& \frac{1}{N}\Theta ^{\left( N\right) }\left( l\right) =\frac{1}{N}\left\{
C-l+\left( N-1\right) \left[ C-\sum\limits_{k=1}^{K+L}kg_{k,j}+\sum%
\limits_{k=K}^{K+L-1}\frac{\left( 1-\beta \right) g_{k,j}}{\left[ 1-\left(
1-\beta \right) g_{k,j}\right] ^{2}}+\frac{g_{K+L,j}}{\left[ 1-g_{K+L,j}%
\right] ^{2}}\right] \right\} \\
& \rightarrow C-\sum\limits_{k=1}^{K+L}kg_{k,j}+\sum\limits_{k=K}^{K+L-1}%
\frac{\left( 1-\beta \right) g_{k,j}}{\left[ 1-\left( 1-\beta \right) g_{k,j}%
\right] ^{2}}+\frac{g_{K+L,j}}{\left[ 1-g_{K+L,j}\right] ^{2}}\overset{\text{%
Def}}{=}\mathbf{\Theta }.
\end{align*}

Let%
\begin{equation*}
\mathbf{A}=\lim_{N\rightarrow\infty}\mathbf{A}_{N},\text{ \ }\mathbf{A}^{%
\text{renting}}=\lim_{N\rightarrow\infty}\mathbf{A}_{N}^{\text{renting}},
\end{equation*}%
\begin{equation*}
\mathbf{A}^{\text{returning}}=\lim_{N\rightarrow\infty}\mathbf{A}_{N}^{\text{%
returning}},\text{ \ }\mathbf{A}^{\text{environment}}=\lim_{N\rightarrow%
\infty}\mathbf{A}_{N}^{\text{environment}}.
\end{equation*}
Then%
\begin{align}
\mathbf{A}f(\mathbf{g})= & -\sum_{j=1}^{m}\lambda_{j}\sum_{k=1}^{K+L}g_{k,j}%
\frac{\partial}{\partial g_{k,j}}f\left( \mathbf{g}\right)
-\alpha\sum_{j=1}^{m}\lambda_{j}\sum_{k=-\left( L-1\right) }^{0}g_{k,j}\frac{%
\partial}{\partial g_{k,j}}f\left( \mathbf{g}\right)  \notag \\
& +\sum_{i=1}^{m}\sum_{j=1}^{m}\sum_{k=-L}^{K+L}g_{k,i}w_{i,j}\left[ -\frac{%
\partial}{\partial g_{k,i}}f(\mathbf{g})+\frac{\partial}{\partial g_{k,j}}f(%
\mathbf{g})\right]  \notag \\
& +\mathbf{\Theta}\left(
\sum_{j=1}^{m}\sum_{l=-L}^{K-1}\mu_{j}g_{l,j}+\beta\sum_{j=1}^{m}%
\sum_{l=K}^{K+L-1}\mu_{j}g_{l,j}\right) \frac{\partial }{\partial g_{k,j}}%
f\left( \mathbf{g}\right) .  \label{oper-2}
\end{align}

Now, we discuss the weak convergence of the sequence $\left\{ \mathbf{Y}%
^{\left( N\right) }\left( t\right) :t\geq0\right\} $ of Markov processes for
$N=1,2,3,\ldots.$ Here, our main purpose is to provide some basic support
for our later study of various convergence involved. To this end, we
consider the random vector $\mathbf{Y}^{\left( N\right) }\left( t\right) $
with samples in $\mathcal{P}\left( \mathbb{D}\left( R_{+},\mathbf{N}\right)
\right) $, where $R_{+}=[0,+\infty)$, $\mathbf{N}=\left( \left( k,j\right)
:-L\leq k\leq K+L,1\leq j\leq m\right) $, $\mathbb{D}\left( R_{+},\mathbf{N}%
\right) $ is the Skorohod space, i.e., the set of mappings which are right
continuous with left-hand limits (in short, C\`{a}dl\`{a}g), and $\mathcal{P}%
\left( \mathbb{\cdot}\right) $ is the set of probability measures defined in
$\mathbb{D}\left( R_{+},\mathbf{N}\right) $. Notice that the convergence in
the Skorohod topology means the convergence in distribution (or weak
convergence) for the Skorohod topology on the space of trajectories. When
the sequence $\left\{ \mathbf{Y}^{(N)}(t),t\geq0\right\} $ of Markov
processes converges in probability (or converges weakly), for the Skorohod
topology, to a given probability vector $\mathbf{Y}(t)$, we write the weak
convergence as $\mathbf{Y}^{\left( N\right) }\left( t\right) \Longrightarrow%
\mathbf{Y}\left( t\right) $ for $t\geq0$, as $N\longrightarrow\infty$.

If $\mathbf{Y}^{\left( N\right) }\left( t\right) \Longrightarrow \mathbf{Y}%
\left( t\right) $ for $t\geq 0$ as $N\longrightarrow \infty $, then it is
easy to see from (\ref{oper-1}) and (\ref{oper-2}) that the transition
probabilities of the Markov process $\left\{ \mathbf{Y}^{\left( N\right)
}\left( t\right) ,t\geq 0\right\} $ with generating operator $\mathbf{A}_{N}$
uniformly converges on any finite time interval to the transition
probabilities of the limiting Markov process $\left\{ \mathbf{Y}\left(
t\right) ,t\geq 0\right\} $ with generating operator $\mathbf{A}$.

Now we consider the limiting behavior of the sequence $\{\mathbf{Y}%
^{(N)}(t),t\geq0\}$ of Markov processes as $N\rightarrow\infty$. To that
end, we first give a system of limiting mean-field equations (\ref{eq-8}) to
(\ref{eq-8-1}) below.

Set
\begin{equation*}
\mathbf{y}\left( t\right) =\lim\limits_{N\rightarrow\infty}\mathbf{y}%
^{\left( N\right) }\left( t\right)
\end{equation*}
and
\begin{equation*}
\mathbf{V}_{\mathbf{y}\left( t\right) }=\lim\limits_{N\rightarrow\infty }%
\mathbf{V}_{\mathbf{y}^{\left( N\right) }\left( t\right) }.
\end{equation*}
Then it follows from (\ref{eq-7}) and (\ref{eq-7-1}) that%
\begin{equation}
\frac{\text{d}}{\text{d}t}\mathbf{y}\left( t\right) =\mathbf{y}\left(
t\right) \mathbf{V}_{\mathbf{y}\left( t\right) },  \label{eq-8}
\end{equation}%
\begin{equation}
\text{ \ \ }\mathbf{y}\left( t\right) e=1,\text{ }\mathbf{y}\left( 0\right) =%
\mathbf{g\in\Omega.}  \label{eq-8-1}
\end{equation}

Note that the convergence in the Skorohod topology means the convergence in
distribution for the Skorohod topology on the space of trajectories. The
following theorem applies the martingale limit theory to discuss the weak
convergence of the sequence $\left\{ \mathbf{Y}^{(N)}(t),t\geq 0\right\} $
of Markov processes as $N$ tends to infinity.

\begin{The}
\label{The:Mart}If $\mathbf{Y}^{(N)}(0)$ converges weakly to $\mathbf{g}%
\in\Omega$ as $N$ tends to infinity, then the sequence $\left\{ \mathbf{Y}%
^{(N)}(t),t\geq0\right\} $ of Markov processes converges in the Skorohod
topology to a solution $\mathbf{y}\left( t\right) $\ to the system of
limiting mean-field equations (\ref{eq-8}) to (\ref{eq-8-1}).
\end{The}

\textbf{Proof: }From the martingale characterization of the Markov jump
process $\left\{ \mathbf{Y}^{(N)}(t),t\geq 0\right\} $, it follows from
Rogers and Williams \cite{Rog:1987, Rog:1994} that for $-L\leq k\leq K+L$
and $1\leq j\leq m$,\newline
\begin{equation*}
M_{k,j}^{(N)}(t)=Y_{k,j}^{\left( N\right) }\left( t\right) -Y_{k,j}^{\left(
N\right) }\left( 0\right) -\int_{0}^{t}\sum_{\Lambda \in \Omega -\left\{
\mathbf{Y}^{(N)}(t)\right\} }\mathcal{Q}^{\left( N\right) }\left( \mathbf{Y}%
^{(N)}(s),\Lambda \right) \left[ \Lambda _{k,j}-Y_{k,j}^{\left( N\right)
}\left( s\right) \right] \text{d}s
\end{equation*}%
is a martingale with respect to the natural filtration associated to the
Poisson processes involved in the renting and returning processes and to the
Markov process of the Markovian environment, where $\mathcal{Q}^{\left(
N\right) }\left( \mathbf{Y}^{(N)}(s),\Lambda \right) $ is the $Q$-matrix of
the Markov jump process $\left\{ \mathbf{Y}^{(N)}(t),t\geq 0\right\} $ whose
expression is given by means of the state change due to the renting and
returning processes as well as the state transitions of the Markovian
environment.

To express the $Q$-matrix $\mathcal{Q}^{\left( N\right) }\left( \mathbf{Y}%
^{(N)}(s),\Lambda\right) $, we analyze three classes of state transitions as
follows:

(1) When a customer arrives at the tagged station to rent a bike, the state
transition rate is given by
\begin{equation*}
q_{k,j;k-1,j}=\left\{
\begin{array}{l}
\lambda_{j},\text{ \ \ }1\leq k\leq K+L,\text{ \ \ \ \ \ }1\leq j\leq m, \\
\lambda_{j}\alpha,\text{ }-\left( L-1\right) \leq k\leq0,\text{ \ }1\leq
j\leq m.%
\end{array}
\right.
\end{equation*}

(2) When a customer returns his bike to the tagged station, the state
transition rate is given by%
\begin{equation*}
q_{k,j;k+1,j}^{\left( N\right) }\left( t\right) =\left\{
\begin{array}{ll}
\frac{\mu _{j}}{N}\left\{ C+\left( N-1\right) \left[ C-\sum%
\limits_{k=1}^{K+L}ky_{k,j}^{\left( N\right) }\left( t\right) \right. \right.
&  \\
\left. \left. +\sum\limits_{k=K}^{K+L-1}\frac{\left( 1-\beta \right)
y_{k,j}^{\left( N\right) }\left( t\right) }{\left[ 1-\left( 1-\beta \right)
y_{k,j}^{\left( N\right) }\left( t\right) \right] ^{2}}+\frac{%
y_{K+L,j}^{\left( N\right) }\left( t\right) }{\left[ 1-y_{K+L,j}^{\left(
N\right) }\left( t\right) \right] ^{2}}\right] \right\} , & \text{\ }-L\leq
l\leq 0, \\
\frac{\mu _{j}}{N}\left\{ C-l+\left( N-1\right) \left[ C-\sum%
\limits_{k=1}^{K+L}ky_{k,j}^{\left( N\right) }\left( t\right) \right. \right.
&  \\
\left. \left. +\sum\limits_{k=1}^{K+L-1}\frac{\left( 1-\beta \right)
y_{k,j}^{\left( N\right) }\left( t\right) }{\left[ 1-\left( 1-\beta \right)
y_{k,j}^{\left( N\right) }\left( t\right) \right] ^{2}}+\frac{%
y_{K+L,j}^{\left( N\right) }\left( t\right) }{\left[ 1-y_{K+L,j}^{\left(
N\right) }\left( t\right) \right] ^{2}}\right] \right\} , & 1\leq l\leq C-1,
\\
\frac{\mu _{j}}{N}\left\{ \left( N-1\right) \left[ C-\sum%
\limits_{k=1}^{K+L}ky_{k,j}^{\left( N\right) }\left( t\right) \right. \right.
&  \\
\left. \left. +\sum\limits_{k=1}^{K+L-1}\frac{\left( 1-\beta \right)
y_{k,j}^{\left( N\right) }\left( t\right) }{\left[ 1-\left( 1-\beta \right)
y_{k,j}^{\left( N\right) }\left( t\right) \right] ^{2}}+\frac{%
y_{K+L,j}^{\left( N\right) }\left( t\right) }{\left[ 1-y_{K+L,j}^{\left(
N\right) }\left( t\right) \right] ^{2}}\right] \right\} , & C\leq l\leq K-1,
\\
\beta \frac{\mu _{j}}{N}\left\{ \left( N-1\right) \left[ C-\sum%
\limits_{k=1}^{K+L}ky_{k,j}^{\left( N\right) }\left( t\right) \right. \right.
&  \\
\left. \left. +\sum\limits_{k=1}^{K+L-1}\frac{\left( 1-\beta \right)
y_{k,j}^{\left( N\right) }\left( t\right) }{\left[ 1-\left( 1-\beta \right)
y_{k,j}^{\left( N\right) }\left( t\right) \right] ^{2}}+\frac{%
y_{K+L,j}^{\left( N\right) }\left( t\right) }{\left[ 1-y_{K+L,j}^{\left(
N\right) }\left( t\right) \right] ^{2}}\right] \right\} , & K\leq l\leq
K+L-1.%
\end{array}%
\right.
\end{equation*}

(3) When the Markovian environment changes from State $i$ to State $j$, the
state transition rate is given by
\begin{equation*}
q_{k,i;k,j}=w_{i,j},\text{ \ }-L\leq k\leq K+L,i\neq j,1\leq i,j\leq m.
\end{equation*}

Based on the above three cases, the $Q$-matrix $\mathcal{Q}^{\left( N\right)
}\left( \mathbf{Y}^{(N)}(s),\Lambda \right) $ is given by%
\begin{align*}
Y_{k,j}^{\left( N\right) }\left( t\right) =&
M_{k,j}^{(N)}(t)+Y_{k,j}^{\left( N\right) }\left( 0\right)
+q_{k+1,j;k,j}\int_{0}^{t}Y_{k+1,j}^{\left( N\right) }\left( s\right) \text{d%
}s \\
& +\int_{0}^{t}q_{k-1,j;k,j}^{\left( N\right) }\left( s\right)
Y_{k-1,j}^{\left( N\right) }\left( s\right) \text{d}s+\sum_{i\neq
j}^{m}w_{i,j}\int_{0}^{t}Y_{k,i}^{\left( N\right) }\left( s\right) \text{d}s
\end{align*}%
Using a similar method to Darling and Norris \cite{Dar:2005, Dar:2008}, it
is easy to see that if $\mathbf{Y}^{(N)}(0)$ converges weakly to $\mathbf{g}%
\in \Omega $ as $N$ tends to infinity, then the sequence $\left\{ \mathbf{Y}%
^{(N)}(t),t\geq 0\right\} $ of Markov processes is tight for the Skorohod
topology, and any limit $\mathbf{Y}(t)$ of $\left\{ \mathbf{Y}%
^{(N)}(t),t\geq 0\right\} $ asymptotically approaches to a single trajectory
identified by a solution $\mathbf{y}\left( t\right) $\ to the system of
limiting mean-field equations (\ref{eq-8}) to (\ref{eq-8-1}). This completes
the proof. \textbf{{\rule{0.08in}{0.08in}}}

\section{A Nonlinear QBD Process}

In this section, we discuss the fixed point of the block-structured system
of limiting mean-field equations (\ref{eq-8}) to (\ref{eq-8-1}), and provide
a mean-field matrix-analytic method which can be used to numerically compute
the fixed point. Furthermore, we study the limiting interchangeability of $%
\mathbf{y}^{(N)}(t)$ as $N\rightarrow \infty $ and $t\rightarrow +\infty $,
that is, the asymptotic independency, which plays a key role in approximate
computation for performance measures of this bike sharing system.

We rewrite the system of limiting mean-field equations (\ref{eq-8}) to (\ref%
{eq-8-1}) as%
\begin{equation*}
\frac{\text{d}}{\text{d}t}\mathbf{y}\left( t\right) =\mathbf{y}\left(
t\right) \mathbf{V}_{\mathbf{y}\left( t\right) }
\end{equation*}%
and%
\begin{equation*}
\text{ \ \ }\mathbf{y}\left( t\right) e=1,\text{ }\mathbf{y}\left( 0\right) =%
\mathbf{g\in \Omega .}
\end{equation*}%
A point $\mathbf{\pi }\in \Omega $ is said to be a fixed point if $%
\lim_{t\rightarrow +\infty }\left[ \frac{\text{d}}{\text{d}t}\mathbf{y}%
\left( t\right) \right] =0$, or%
\begin{equation*}
\left[ \mathbf{y}\left( t\right) \mathbf{V}_{\mathbf{y}\left( t\right) }%
\right] _{|\mathbf{y}\left( t\right) =\mathbf{\pi }}=0.
\end{equation*}%
Thus, we have%
\begin{equation}
\mathbf{\pi V}_{\pi }=0  \label{FixedP-1}
\end{equation}%
and%
\begin{equation}
\mathbf{\pi }e=1.  \label{FixedP-2}
\end{equation}

Now, we provide a mean-field matrix-analytic method to compute the fixed
point $\mathbf{\pi }$ from the system of nonlinear equations: $\mathbf{\pi V}%
_{\pi }=0$ and\ $\mathbf{\pi e}=1$. To this end, it is necessary to explore
the block structure of the system of nonlinear equations. Hence this gives a
nonlinear QBD process so that the $RG$-factorizations given by Li \cite%
{Li:2010} are applicable in our later analysis.

Let%
\begin{equation*}
\xi _{k,j}=\lim\limits_{t\rightarrow +\infty }\lim\limits_{N\rightarrow
\infty }\xi _{k,j}^{\left( N\right) }\left( t\right) ,\text{ \ }-L\leq k\leq
K+L,1\leq j\leq m.
\end{equation*}%
Then%
\begin{equation*}
\xi _{k,j}=\left\{
\begin{array}{l}
\mu _{j}\zeta _{j},\text{ }-L\leq k\leq K-1,\text{ \ \ \ \ }1\leq j\leq m,
\\
\beta \mu _{j}\zeta _{j},\text{\ \ }K\leq k\leq K+L-1,\text{\ }1\leq j\leq m,%
\end{array}%
\right.
\end{equation*}%
where
\begin{equation*}
\zeta _{j}=C-\sum\limits_{k=1}^{K+L}k\pi _{k,j}+\sum\limits_{k=K}^{K+L-1}%
\frac{\left( 1-\beta \right) \pi _{k,j}}{\left[ 1-\left( 1-\beta \right) \pi
_{k,j}\right] ^{2}}+\frac{\pi _{K+L,j}}{\left[ 1-\pi _{K+L,j}\right] ^{2}}.
\end{equation*}%
Thus for $K\leq k\leq K+L-1$,%
\begin{equation*}
\Psi _{k}=\left(
\begin{array}{cccc}
0 & \beta \mu _{1}\zeta _{1}w_{1,m} & \cdots & \beta \mu _{1}\zeta
_{1}w_{1,m} \\
\beta \mu _{2}\zeta _{2}w_{2,1} & 0 & \cdots & \beta \mu _{2}\zeta
_{2}w_{2,m} \\
\vdots & \vdots & \ddots & \vdots \\
\beta \mu _{m}\zeta _{m}w_{m,1} & \beta \mu _{m}\zeta _{m}w_{m,2} & \cdots &
0%
\end{array}%
\right) \overset{\text{Def}}{=}\Psi \left( \beta \right) ,
\end{equation*}%
\begin{equation*}
\widehat{\Psi }_{k}=\text{diag}\left( \beta \mu _{1}\zeta _{1}w_{1,1},\beta
\mu _{2}\zeta _{2}w_{2,2},\ldots ,\beta \mu _{m}\zeta _{m}w_{m,m}\right)
\overset{\text{Def}}{=}\widehat{\Psi }\left( \beta \right) ;
\end{equation*}%
and for $-L\leq k\leq K-1,$%
\begin{equation*}
\Psi _{k}=\Psi \left( 1\right)
\end{equation*}

and%
\begin{equation*}
\widehat{\Psi}_{k}=\widehat{\Psi}\left( 1\right) .
\end{equation*}

Based on the above analysis, we can summarize the special structures of the
four functions $\Psi _{k}$,$\widehat{\Psi }_{k}$,$\Phi _{k}$ and $\widehat{%
\Phi }_{k}$ in Table 1.
\begin{table}[tbp]
\caption{The special structures of the four functions $\Psi _{k}$,$\widehat{%
\Psi }_{k}$,$\Phi _{k}$ and $\widehat{\Phi }_{k}$}\centering
\begin{tabular}{|c|c|c|c|c|c|}
\hline
\thead{$K$} & \thead{ $-L$ } & \thead{$-L+1\leq k\leq 0$ } & \thead{$1\leq
k\leq K-1$ } & \thead{$K\leq k\leq K+L-1$} & \thead{$K+L$} \\ \hline
$\Psi _{k}$ & $\Psi \left( 1\right) $ & $\Psi \left( 1\right) $ & $\Psi
\left( 1\right) $ & $\Psi \left( \beta \right) $ & null \\ \hline
$\widehat{\Psi }_{k}$ & $\widehat{\Psi }\left( 1\right) $ & $\widehat{\Psi }%
\left( 1\right) $ & $\widehat{\Psi }\left( 1\right) $ & $\widehat{\Psi }%
\left( \beta \right) $ & null \\ \hline
$\Phi _{k}$ & null & $\Phi\left( \alpha \right) $ & $\Phi \left( 1\right) $
& $\Phi \left( 1\right) $ & $\Phi \left( 1\right) $ \\ \hline
$\widehat{\Phi }_{k}$ & null & $\widehat{\Phi }\left( \alpha \right) $ & $%
\widehat{\Phi }\left( 1\right) $ & $\widehat{\Phi }\left( 1\right) $ & $%
\widehat{\Phi }\left( 1\right) $ \\ \hline
\end{tabular}
\end{table}

By observing Table 1, it is easy to check from (\ref{eq-7-2}) that as $%
N\rightarrow \infty $ and $t\rightarrow +\infty $%
\begin{equation}
\mathbf{V}_{\pi }=\left(
\begin{array}{ccc}
B_{1,1} & B_{1,2} &  \\
B_{2,1} & B_{2,2} & B_{2,3} \\
& B_{3,2} & B_{3,3}%
\end{array}%
\right) ,  \label{eq-10}
\end{equation}%
where%
\begin{equation*}
B_{1,1}=\left(
\begin{array}{ccccc}
\widehat{\Psi }_{-L} & \Psi _{-L} &  &  &  \\
\Phi _{-\left( L-1\right) } & \widehat{\Phi }_{-\left( L-1\right) }+\widehat{%
\Psi }_{-\left( L-1\right) } & \Psi _{-\left( L-1\right) } &  &  \\
& \ddots & \ddots & \ddots &  \\
&  & \Phi _{-1} & \widehat{\Phi }_{-1}+\widehat{\Psi }_{-1} & \Psi _{-1} \\
&  &  & \Phi _{0} & \widehat{\Phi }_{0}+\widehat{\Psi }_{0}%
\end{array}%
\right) ,
\end{equation*}%
\begin{equation*}
B_{1,2}=\left(
\begin{array}{cccc}
&  &  &  \\
&  &  &  \\
&  &  &  \\
\Psi _{0} &  &  &
\end{array}%
\right) ,\text{ \ }B_{2,1}=\left(
\begin{array}{cccc}
&  &  & \Phi _{1} \\
&  &  &  \\
&  &  &  \\
&  &  &
\end{array}%
\right) ,
\end{equation*}%
\begin{equation*}
B_{2,2}=\left(
\begin{array}{ccccc}
\widehat{\Phi }_{1}+\widehat{\Psi }_{1} & \Psi _{1} &  &  &  \\
\Phi _{2} & \widehat{\Phi }_{2}+\widehat{\Psi }_{2} & \Psi _{2} &  &  \\
& \ddots & \ddots & \ddots &  \\
&  & \Phi _{K-2} & \widehat{\Phi }_{K-2}+\widehat{\Psi }_{K-2} & \Psi _{K-2}
\\
&  &  & \Phi _{K-1} & \widehat{\Phi }_{K-1}+\widehat{\Psi }_{K-1}%
\end{array}%
\right) ,
\end{equation*}%
\begin{equation*}
B_{2,3}=\left(
\begin{array}{cccc}
&  &  &  \\
&  &  &  \\
&  &  &  \\
\Psi _{K-1} &  &  &
\end{array}%
\right) ,\text{ \ }B_{3,2}=\left(
\begin{array}{cccc}
&  &  & \Phi _{K} \\
&  &  &  \\
&  &  &  \\
&  &  &
\end{array}%
\right) ,
\end{equation*}%
\begin{equation*}
B_{3,3}=\left(
\begin{array}{ccccc}
\widehat{\Phi }_{K}+\widehat{\Psi }_{K} & \Psi _{K} &  &  &  \\
\Phi _{K+1} & \widehat{\Phi }_{K+1}+\widehat{\Psi }_{K+1} & \Psi _{K+1} &  &
\\
& \ddots & \ddots & \ddots &  \\
&  & \Phi _{K+L-1} & \widehat{\Phi }_{K+L-1}+\widehat{\Psi }_{K+L-1} & \Psi
_{K+L-1} \\
&  &  & \Phi _{K+L} & \widehat{\Phi }_{K+L}%
\end{array}%
\right) .
\end{equation*}%
Thus it follows from (\ref{eq-10}) that

\begin{equation*}
\mathbf{V}_{\pi}=\left(
\begin{array}{ccccc}
\widehat{\Psi}_{-L} & \Psi_{-L} &  &  &  \\
\Phi_{-\left( L-1\right) } & \widehat{\Phi}_{-\left( L-1\right) }+\widehat{%
\Psi}_{-\left( L-1\right) } & \Psi_{-\left( L-1\right) } &  &  \\
& \ddots & \ddots & \ddots &  \\
&  & \Phi_{K+L-1} & \widehat{\Phi}_{K+L-1}+\widehat{\Psi}_{K+L-1} &
\Psi_{K+L-1} \\
&  &  & \Phi_{K+L} & \widehat{\Phi}_{K+L}%
\end{array}
\right) .
\end{equation*}

Based on the nonlinear QBD process $\mathbf{V}_{\pi }$, we write
\begin{equation*}
\mathbf{\pi }=\left( \mathbf{\pi }_{-L},\mathbf{\pi }_{-L+1},\ldots ,\mathbf{%
\pi }_{K+L-1},\mathbf{\pi }_{K+L}\right) .
\end{equation*}%
Now, we use the LU-type $RG$-factorization given in Subsection 1.3.2 of
Chapter one in Li \cite{Li:2010} (see Pages 25 and 26), and write the
LU-type R-measure as%
\begin{equation*}
\mathbf{R}_{-L+1}\left( \mathbf{\pi }\right) =-\Phi _{-L+1}\left( \widehat{%
\Psi }_{-L}\right) ^{-1},
\end{equation*}%
\begin{equation*}
\mathbf{R}_{-L+2}\left( \mathbf{\pi }\right) =-\Phi _{-L+2}\left[ \mathbf{R}%
_{-L+1}\left( \mathbf{\pi }\right) \Psi _{-L}+\left( \widehat{\Phi }_{-L+1}+%
\widehat{\Psi }_{-L+1}\right) \right] ^{-1},
\end{equation*}%
for $-L+2\leq k\leq K+L$,%
\begin{equation*}
\mathbf{R}_{k}\left( \mathbf{\pi }\right) =-\Phi _{k}\left[ \mathbf{R}%
_{k-1}\left( \mathbf{\pi }\right) \Psi _{k-2}+\left( \widehat{\Phi }_{k-1}+%
\widehat{\Psi }_{k-1}\right) \right] ^{-1}.
\end{equation*}

At the same time, the infinitesimal generator of the censored Markov chain
to level $K+L$ is given by%
\begin{equation*}
\mathbf{\Xi }_{K+L}=\mathbf{R}_{K+L}\left( \mathbf{\pi }\right) \Psi
_{K+L-1}+\widehat{\Phi }_{K+L}.
\end{equation*}%
For the fixed point $\mathbf{\pi }=\left( \mathbf{\pi }_{-L},\mathbf{\pi }%
_{-L+1},\ldots ,\mathbf{\pi }_{K+L-1},\mathbf{\pi }_{K+L}\right) $, it
follows from (1.25) in Subsection 1.3.4.1 of Li \cite{Li:2010} (see Page 30)
that for $k=-L$,%
\begin{equation*}
\mathbf{\pi }_{-L}=\mathbf{\pi }_{-L+1}\mathbf{R}_{-L+1}\left( \mathbf{\pi }%
\right) ,
\end{equation*}%
and for $-L+1\leq k\leq K+L-1$,%
\begin{align*}
\mathbf{\pi }_{k}& =\mathbf{\pi }_{k+1}\mathbf{R}_{k+1}\left( \mathbf{\pi }%
\right) \\
& =\mathbf{\pi }_{K+L}\mathbf{R}_{K+L}\left( \mathbf{\pi }\right) \mathbf{R}%
_{K+L-1}\left( \mathbf{\pi }\right) \mathbf{R}_{K+L-2}\left( \mathbf{\pi }%
\right) \cdots \mathbf{R}_{k+1}\left( \mathbf{\pi }\right) ,
\end{align*}%
where the vector $\mathbf{\pi }_{K+L}$ is a solution to the systems of
nonlinear equations $\mathbf{\pi }_{K+L}\mathbf{\Xi }_{K+L}=0$ and $\mathbf{%
\pi }_{K+L}\left[ I+\sum_{k=-L}^{K+L-1}\mathbf{R}_{K+L}\left( \mathbf{\pi }%
\right) \mathbf{R}_{K+L-1}\left( \mathbf{\pi }\right) \mathbf{R}%
_{K+L-2}\left( \mathbf{\pi }\right) \cdots \mathbf{R}_{k+1}\left( \mathbf{%
\pi }\right) \right] e=1$.

The following theorem is a summarization of the above analysis, and its
proof is easy to only check the system of nonlinear equations $\mathbf{\pi V}%
_{\pi }=0$ and $\mathbf{\pi }e=1$. Thus we omit the proof here.

\begin{The}
\label{The:FixedP}The fixed point $\mathbf{\pi }$ is a solution to the
vector system of nonlinear equations%
\begin{align}
\mathbf{\pi }=& \left( \mathbf{\pi }_{K+L}\mathbf{R}_{K+L}\left( \mathbf{\pi
}\right) \mathbf{R}_{K+L-1}\left( \mathbf{\pi }\right) \mathbf{R}%
_{K+L-2}\left( \mathbf{\pi }\right) \cdots \mathbf{R}_{-L+1}\left( \mathbf{%
\pi }\right) ,\right.  \notag \\
& \left. \mathbf{\pi }_{K+L}\mathbf{R}_{K+L}\left( \mathbf{\pi }\right)
\mathbf{R}_{K+L-1}\left( \mathbf{\pi }\right) \mathbf{R}_{K+L-2}\left(
\mathbf{\pi }\right) \cdots \mathbf{R}_{-L+2}\left( \mathbf{\pi }\right)
,\right.  \notag \\
& \left. \ldots ,\mathbf{\pi }_{K+L}\mathbf{R}_{K+L}\left( \mathbf{\pi }%
\right) \mathbf{R}_{K+L-1}\left( \mathbf{\pi }\right) ,\mathbf{\pi }_{K+L}%
\mathbf{R}_{K+L}\left( \mathbf{\pi }\right) ,\mathbf{\pi }_{K+L}\right) ,
\label{Fix-1}
\end{align}%
\begin{equation}
\mathbf{\pi }_{K+L}\left[ \mathbf{R}_{K+L}\left( \mathbf{\pi }\right) \Psi
_{K+L-1}+\left( \widehat{\Phi }_{K+L}+\widehat{\Psi }_{K+L}\right) \right] =0
\label{Fix-2}
\end{equation}%
and%
\begin{equation}
\mathbf{\pi }_{K+L}\left[ I+\sum_{k=-L}^{K+L-1}\mathbf{R}_{K+L}\left(
\mathbf{\pi }\right) \mathbf{R}_{K+L-1}\left( \mathbf{\pi }\right) \mathbf{R}%
_{K+L-2}\left( \mathbf{\pi }\right) \cdots \mathbf{R}_{k+1}\left( \mathbf{%
\pi }\right) \right] e=1.  \label{Fix-3}
\end{equation}
\end{The}

It is easy to see that although the system of nonlinear equations: $\mathbf{%
\pi V}_{\pi }=0$ and $\mathbf{\pi }e=1$, are equivalent to the vector system
of nonlinear equations (\ref{Fix-1}), (\ref{Fix-2}) and (\ref{Fix-3}),
Theorem \ref{The:FixedP} can be used to design different algorithms to
numerically compute the fixed points $\mathbf{\pi }$. This is indicated in
the next section with six numerical examples. Reader may also refer to Li %
\cite{Li:2016a} and Li \emph{et al}. \cite{Li:2016c} for some nearby
research.

In what follows we discuss the mean-field limit of the empirical measure
process of the bike sharing system as the number $N$ of stations and time $t$%
\ go to infinity, and show that the fixed point is unique from the system of
nonlinear equations: $\mathbf{\pi V}_{\pi }=0$ and $\mathbf{\pi }e=1$. It is
worthwhile to note that the uniqueness of the fixed point guarantees the
\textit{asymptotic independence} of the queueing processes describing the
numbers of bikes at the $N$ stations as $N\rightarrow \infty $, also known
as the \textit{propagation of chaos}.

For the unique fixed point $\mathbf{\pi }$, we discuss the limiting
interchangeability of the probability vector $\mathbf{y}^{\left( N\right)
}\left( t,\mathbf{g}\right) $ as $N\rightarrow \infty $ and $t\rightarrow
+\infty $, where $\mathbf{y}^{\left( N\right) }\left( 0,\mathbf{g}\right) =%
\mathbf{g\in \Omega }$. Note that the limiting interchangeability is
necessary in many practical applications when using the stationary
probabilities (that is, the fixed point) of the limiting process to give an
effective approximation for performance analysis of this bike sharing system.

The following theorem gives the limit of the vector $\mathbf{y}(t,\mathbf{g}%
) $ as $t\rightarrow +\infty $, that is,%
\begin{equation*}
\mathbf{y}(t,\mathbf{g})=\lim_{N\rightarrow \infty }\mathbf{y}^{\left(
N\right) }(t,\mathbf{g})
\end{equation*}%
and%
\begin{equation*}
\lim_{t\rightarrow +\infty }\mathbf{y}(t,\mathbf{g})=\lim_{t\rightarrow
+\infty }\lim_{N\rightarrow \infty }\mathbf{y}^{\left( N\right) }(t,\mathbf{g%
}).
\end{equation*}

\begin{The}
\label{The:Limit1}For any $\mathbf{g}\in\Omega$%
\begin{equation*}
\lim_{t\rightarrow+\infty}\mathbf{y}(t,\mathbf{g})=\mathbf{\pi}.
\end{equation*}
Furthermore, there exists a unique probability measure $\varphi$ on $\Omega$%
, which is invariant under the map $\mathbf{g}\longmapsto\mathbf{y}(t,%
\mathbf{g})$, that is, for any continuous function $f:$ $\Omega \rightarrow%
\mathbf{R}$ and $t>0$%
\begin{equation*}
\int_{\Omega}f(\mathbf{g})\text{d}\varphi(\mathbf{g})=\int_{\Omega }f(%
\mathbf{y}(t,\mathbf{g}))\text{d}\varphi(\mathbf{g}).
\end{equation*}
Also, $\varphi=\delta_{\mathbf{\pi}}$ is the probability measure
concentrated at the fixed point $\mathbf{\pi}$.
\end{The}

\textbf{Proof:} It is seen from Theorem \ref{The:Mart} that as $t\rightarrow
+\infty$,\ the limit of $\mathbf{y}(t,\mathbf{g})$ exists on $\Omega$, and
it is also a solution on $\Omega$ to the system of nonlinear equations (\ref%
{FixedP-1}) and (\ref{FixedP-2}). Since $\mathbf{y}(t,\mathbf{g})$ is the
unique solution to the system of limiting mean-field equations (\ref{eq-8})
and (\ref{eq-8-1}), the vector $\lim_{t\rightarrow+\infty}\mathbf{y}(t,%
\mathbf{g})$ is also a solution to the system of nonlinear equations (\ref%
{FixedP-1}) and (\ref{FixedP-2}). Note that $\mathbf{\pi}$ is the unique
solution to the system of nonlinear equations (\ref{FixedP-1}) and (\ref%
{FixedP-2}), hence we obtain that $\lim_{t\rightarrow+\infty}\mathbf{y}(t,%
\mathbf{g})=\mathbf{\pi}$. The second statement in this theorem can be
immediately given by the probability measure of the limiting process $%
\left\{ \mathbf{Y}(t),t\geq0\right\} $ on state space $\Omega$. This
completes the proof. \textbf{{\rule{0.08in}{0.08in}}}

The following theorem indicates the weak convergence of the sequence $%
\left\{ \varphi_{N}\right\} $ of stationary probability distributions for
the sequence $\left\{ \mathbf{Y}^{(N)}(t),t\geq0\right\} $ of Markov
processes to the probability measure concentrated at the fixed point $%
\mathbf{\pi}$.

\begin{The}
\label{The:Limit2}(1) For a fixed number $N=1,2,3,\ldots$, the Markov
process $\left\{ \mathbf{Y}^{(N)}(t),t\geq0\right\} $ is positive recurrent,
and has a unique invariant distribution $\varphi_{N}$.

(2) $\left\{ \varphi_{N}\right\} $ weakly converges to $\delta_{\pi}$, that
is, for any continuous function $f:$ $\Omega\rightarrow\mathbf{R}$%
\begin{equation*}
\lim_{N\rightarrow\infty}E_{\varphi_{N}}\left[ f(\mathbf{g})\right] =f\left(
\mathbf{\pi}\right) .
\end{equation*}
\end{The}

\textbf{Proof:} (1) From Theorem 3, this bike sharing system of $N$
identical stations is stable, hence this bike sharing system has a unique
invariant distribution $\varphi_{N}$.

(2) Since $\Omega$ is compact under the metric $\rho\left( \mathbf{g},%
\mathbf{g}^{\prime}\right) $, so it is the set $\mathcal{P}\left(
\Omega\right) $ of probability measures. Hence the sequence $\left\{
\varphi_{N}\right\} $ of invariant distributions has limiting points. A
similar analysis to the proof of Theorem 5 in Martin and Suhov \cite%
{Mar:1999} shows that $\left\{ \varphi_{N}\right\} $ weakly converges to $%
\delta _{\mathbf{\pi}}$ and $\lim_{N\rightarrow\infty}E_{\varphi_{N}}\left[
f(\mathbf{g})\right] =f\left( \mathbf{\pi}\right) $. This completes the
proof. \textbf{{\rule{0.08in}{0.08in}}}

Based on Theorems \ref{The:Limit1} and \ref{The:Limit2}, we obtain a useful
relation as follows%
\begin{equation*}
\lim_{t\rightarrow +\infty }\lim_{N\rightarrow \infty }\mathbf{y}^{\left(
N\right) }(t,\mathbf{g})=\lim_{N\rightarrow \infty }\lim_{t\rightarrow
+\infty }\mathbf{y}^{\left( N\right) }(t,\mathbf{g})=\mathbf{\pi }.
\end{equation*}%
Therefore, we have%
\begin{equation*}
\lim_{\substack{ N\rightarrow \infty  \\ t\rightarrow +\infty }}\mathbf{y}%
^{\left( N\right) }(t,\mathbf{g})=\mathbf{\pi },
\end{equation*}%
which justifies the exchangeability of the limits of $N\rightarrow \infty $
and $t\rightarrow +\infty $.

Finally, we further show the asymptotic independence (or propagation of
chaos) of the queueing processes of this bike sharing system for each $%
k=2,3,\ldots ,N$\ as follows:%
\begin{align*}
& \lim_{t\rightarrow +\infty }\lim_{N\rightarrow \infty }P\left\{
X_{1}^{\left( N\right) }\left( t\right) =n_{1},J_{1}\left( t\right)
=j_{1};\ldots ;X_{k}^{\left( N\right) }\left( t\right) =n_{k},J_{k}\left(
t\right) =j_{k}\right\} \\
& =\lim_{N\rightarrow \infty }\lim_{t\rightarrow +\infty }P\left\{
X_{1}^{\left( N\right) }\left( t\right) =n_{1},J_{1}\left( t\right)
=j_{1};\ldots ;X_{k}^{\left( N\right) }\left( t\right) =n_{k},J_{k}\left(
t\right) =j_{k}\right\} \\
& =\prod_{l=1}^{k}\pi _{n_{l},j_{l}}
\end{align*}%
and%
\begin{align*}
& \lim_{N\rightarrow \infty }\lim_{t\rightarrow +\infty }\frac{1}{t}%
\int_{0}^{t}\mathbf{1}_{\left\{ X_{1}^{\left( N\right) }\left( t\right)
=n_{1},J_{1}\left( t\right) =j_{1};\ldots ;X_{k}^{\left( N\right) }\left(
t\right) =n_{k},J_{k}\left( t\right) =j_{k}\right\} }\text{d}t \\
& =\lim_{t\rightarrow +\infty }\lim_{N\rightarrow \infty }\frac{1}{t}%
\int_{0}^{t}\mathbf{1}_{\left\{ X_{1}^{\left( N\right) }\left( t\right)
=n_{1},J_{1}\left( t\right) =j_{1};\ldots ;X_{k}^{\left( N\right) }\left(
t\right) =n_{k},J_{k}\left( t\right) =j_{k}\right\} }\text{d}t \\
& =\prod_{l=1}^{k}\pi _{n_{l},j_{l}}\text{ \ \ \ \ \ \ \ \ a.s. }
\end{align*}%
It is obvious that the asymptotic independence needs to hold for each subset
of the $N$ same stations. Based on this, it is easy to see that the two
types of limits may be used as an approximate computation for performance
measures of this bike sharing system, hence this demonstrates the key role
played by the asymptotic independence.

\section{Numerical Analysis}

In this section, we first use the fixed point to express interesting
performance measures of this bike sharing system, such as, the stationary
average number of bikes at the tagged station, the stationary
strong-probability of problematic stations, the stationary weak-probability
of problematic stations, and impact of the user's finite waiting rooms on
system performance. Then we use six numerical examples to demonstrate how
the performance measures depend on some key parameters of this bike sharing
system. Therefore, this paper provides numerical solution in the study of
more general bike sharing systems by means of the nonlinear QBD processes.

\subsection{Performance measures}

Using the fixed point $\mathbf{\pi}=\left( \pi_{-L},\pi_{-L+1},\ldots,\pi
_{0},\pi_{1},\ldots,\pi_{K+L-1},\pi_{K+L}\right) $ where $\pi_{k}=\left(
\pi_{k,1},\pi_{k,2},\right. $ $\left. \ldots,\pi_{k,m}\right) $ and $%
\pi_{k}e=\sum_{j=1}^{m}\pi_{k,j}$, we provide some interesting performance
measures of this bike sharing system from a practical point of view as
follows:

\textbf{(1) The stationary average number of bikes parked at the tagged
station}%
\begin{equation*}
E\left[ Q\right] =\sum_{k=1}^{K+L}k\pi _{k}e.
\end{equation*}

\textbf{(2-1) The stationary average number of waiting places used by
customers who are renting bikes at the tagged station}%
\begin{equation*}
E\left[ N_{1}\right] =\sum_{k=-L}^{-1}\left( -k\right) \pi _{k}e.
\end{equation*}

\textbf{(2-2) The stationary average number of waiting places used by
customers who are returning bikes at the tagged station}%
\begin{equation*}
E\left[ N_{2}\right] =\sum_{k=K+1}^{K+L}\left( k-K\right) \pi _{k}e.
\end{equation*}

\textbf{(2-3) The maximal stationary average number of waiting places used
at tagged station}%
\begin{equation*}
E\left[ N\right] =\max \left\{ E\left[ N_{1}\right] ,E\left[ N_{2}\right]
\right\} =\max \left\{ \sum_{k=-L}^{-1}\left( -k\right) \pi
_{k}e,\sum_{k=K+1}^{K+L}\left( k-K\right) \pi _{k}e\right\} .
\end{equation*}%
Since $E\left[ N_{1}\right] $ and $E\left[ N_{2}\right] $ can not exist
simultaneously, $E\left[ N\right] $ is useful for synthetically designing
the user's finite waiting rooms of this bike sharing system.

\textbf{(3) The stationary strong-probability of problematic stations}

The strong-probability of problematic stations is a probability that either
there is both no bike ($-L\leq k\leq 0$) and no empty waiting place ($k=-L$)
when renting a bike, or there is both no parking place ($K\leq k\leq K+L$)
and no empty waiting place ($k=K+L$) when returning a bike. Thus the
stationary strong-probability of problematic stations is given by
\begin{equation*}
\mathbf{p}_{s}=\pi _{-L}e+\pi _{K+L}e.
\end{equation*}

We consider the effect of the size of waiting places on the
strong-probability of problematic stations. Let%
\begin{equation*}
\upsilon =\frac{\pi _{-L}e+\pi _{K+L}e}{\left( \pi _{-L}e+\pi _{K+L}e\right)
_{|L=0}}.
\end{equation*}%
Then $\upsilon $ denotes the improved efficiency of problematic stations due
to introduction of the user's waiting room of size $L>0$.

\textbf{(4) The stationary weak-probability of problematic stations}

The weak-probability of problematic stations is the probability that either
there is no bike ($-L\leq k\leq0$) when renting a bike, or there is no
parking place ($K\leq k\leq K+L$) when returning a bike. Thus the stationary
weak-probability of problematic stations is given by%
\begin{equation*}
\mathbf{p}_{w}=\sum_{k=-L}^{0}\pi_{k}e+\sum_{k=K}^{K+L}\pi_{k}e.
\end{equation*}

\subsection{Numerical examples}

Now, we use six numerical examples to show how the performance measures
depend on some key parameters of this bike sharing system. As illustrated in
the following figures, the Markovian environment motivates us to propose the
mean-field matrix-analytic method which is necessarily developed as some
effective numerical solution in the study of bike sharing systems.

In the following examples one to four, we take some common parameters as
follows:%
\begin{equation*}
K=20,C=10,L=5,\lambda _{2}=50,\mu _{2}=20,\alpha =0.5,\beta
=0.5,m=2,w=\left(
\begin{array}{cc}
-1 & 1 \\
1 & -1%
\end{array}%
\right) ;
\end{equation*}%
while the other parameters are conceretely chosen in each example for the
target of specific observation.

\textbf{Example one: Analysis of }$E\left[ Q\right] $

The left of figure 5 indicates how the stationary average number\textbf{\ }$E%
\left[ Q\right] $ of bikes at the tagged station depends on $\lambda _{1}\in
\left( 30,45\right) $ when $\mu _{1}=25$, $30$, $35$ and $40$, respectively.
It is seen that $E\left[ Q\right] $ decreases as $\lambda _{1}$ increases
but it increases as $\mu _{1}$ increases. The right of figure 5 shows how $E%
\left[ Q\right] $ depends on $\mu _{1}\in \left( 25,40\right) $ when $%
\lambda _{1}=25$, $30$, $35$ and $40$, respectively. It is seen that $E\left[
Q\right] $ increases as $\mu _{1}$ increases but it decreases as $\lambda
_{1}$ increases. Note that these numerical results may intuitively be
understood as follows: The number of rented bikes increases as $\lambda _{1}$
increases, thus $E\left[ Q\right] $ decreases; while the number of returned
bikes increases as $\mu _{1}$ increases, this shows that $E\left[ Q\right] $
increases as $\mu _{1}$ increases.

\begin{figure}[ptb]
\setlength{\abovecaptionskip}{0.cm} \setlength{\belowcaptionskip}{-0.cm} %
\centering              \includegraphics[width=7cm]{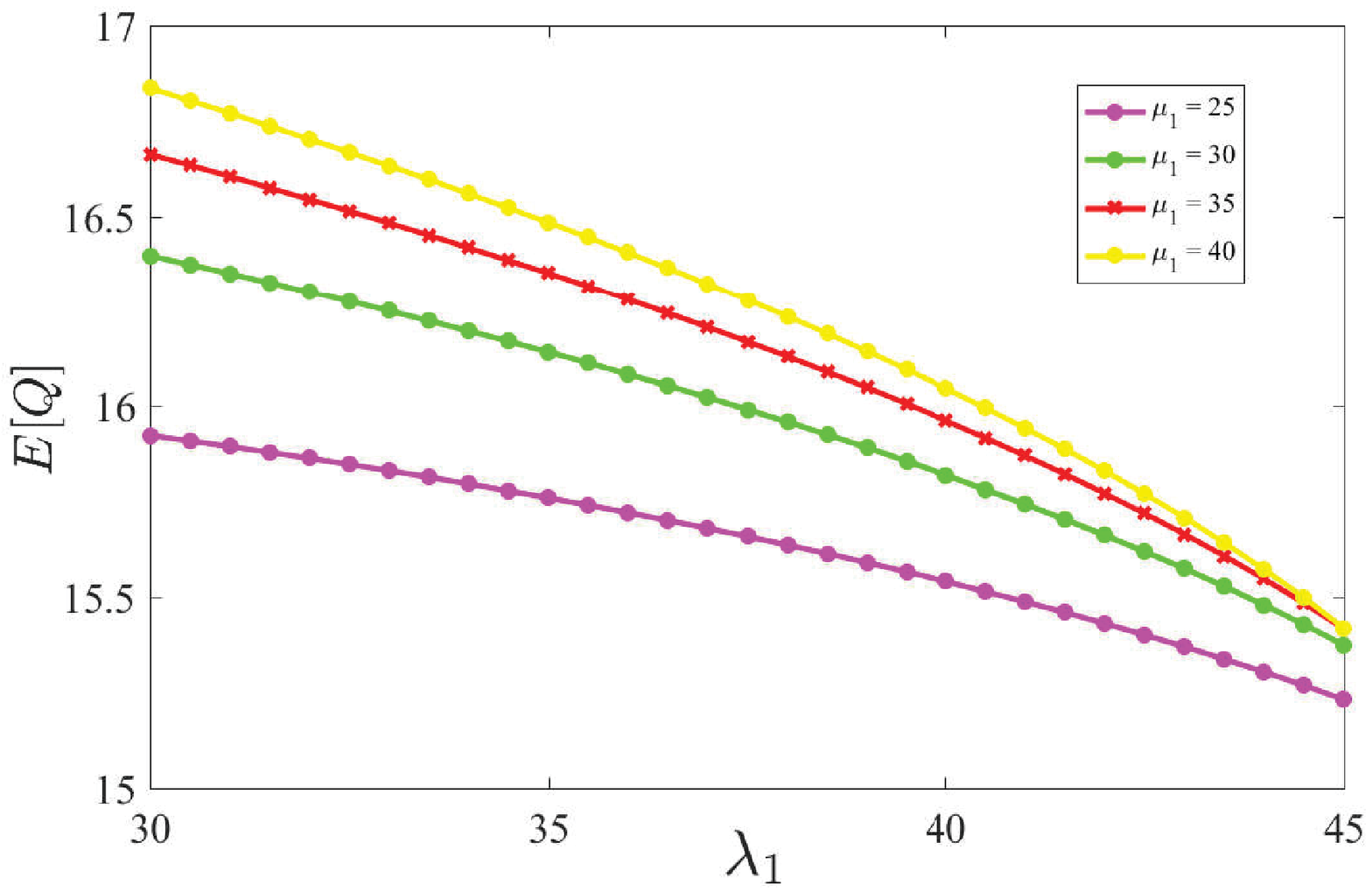} %
\includegraphics[width=7cm]{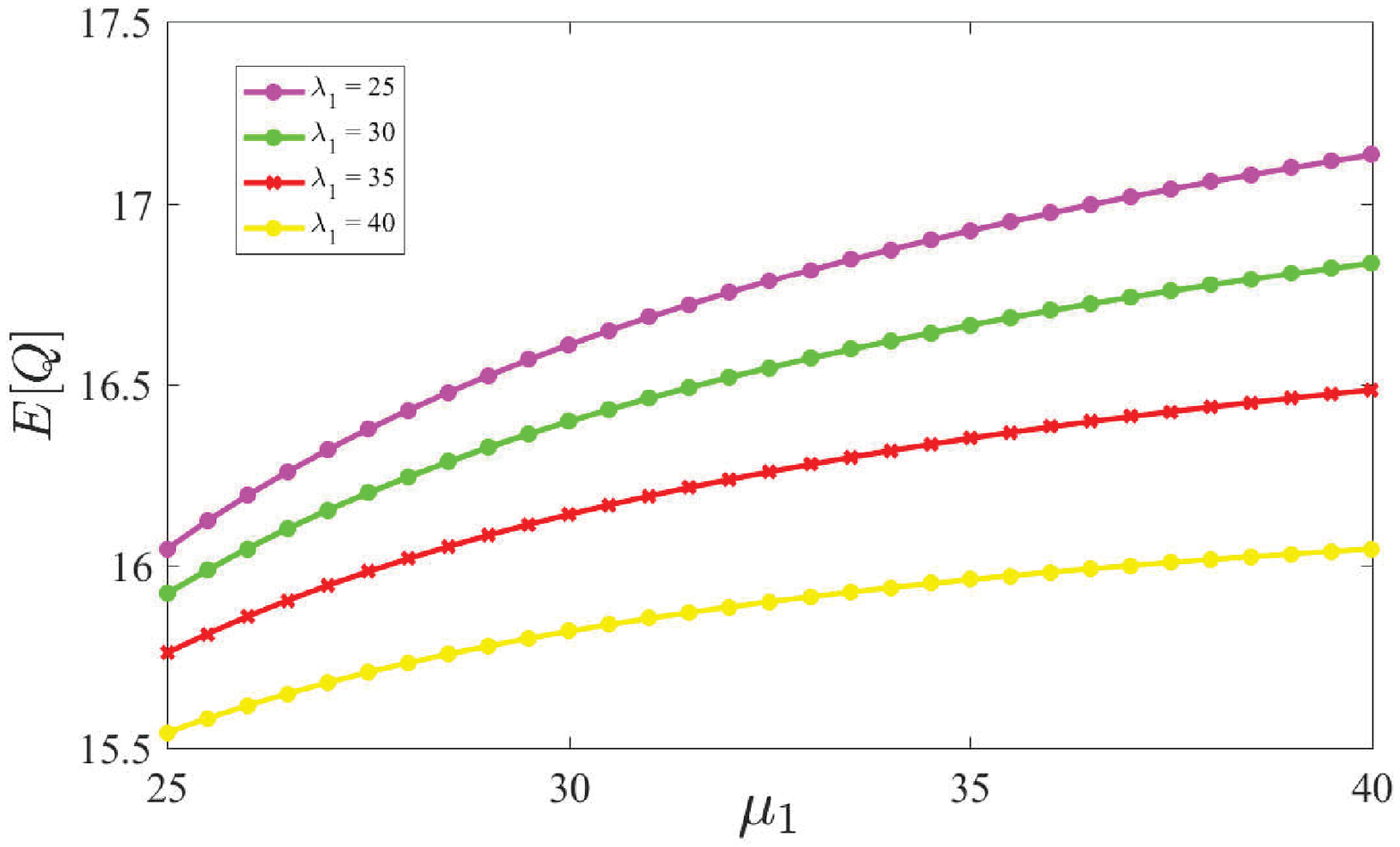} \newline
\caption{$E\left[ Q\right] $ vs. $\protect\lambda_{1}$ and $\protect\mu_{1}$}
\label{figure:figure-5}
\end{figure}

\textbf{Example two: Analysis of }$E\left[ N\right] $

The left of figure 6 shows how the maximal stationary average number $E\left[
N\right] $ of waiting customers at the tagged station depends on $\lambda
_{1}\in \left( 30,45\right) $ when $\mu _{1}=25$, $30$, $35$ and $40$,
respectively. It is seen that $E\left[ N\right] $ decreases as $\lambda _{1}$
increases but it increases as $\mu _{1}$ increases. The right of figure 6
shows how $E\left[ N\right] $ depends on $\mu _{1}\in \left( 25,40\right) $
when $\lambda _{1}=25$, $30$, $35$ and $40$, respectively. It is seen that $E%
\left[ N\right] $ increases as $\mu _{1}$ increases but it decreases as $%
\lambda _{1}$ increases. Note that the number of rented bikes increases as $%
\lambda _{1}$ increases, thus $E\left[ N_{1}\right] $ increases but $E\left[
N_{2}\right] $ decreases. On the other hand, the number of returned bikes
increases as $\mu _{1}$ increases, so $E\left[ N_{2}\right] $ increases but $%
E\left[ N_{1}\right] $ decreases. Based on this, it is clear that $E\left[
N_{2}\right] $ has more impact on $E\left[ N\right] $ than $E\left[ N_{1}%
\right] $ under the present parameter design.

\begin{figure}[ptb]
\setlength{\abovecaptionskip}{0.cm} \setlength{\belowcaptionskip}{-0.cm} %
\centering              \includegraphics[width=7cm]{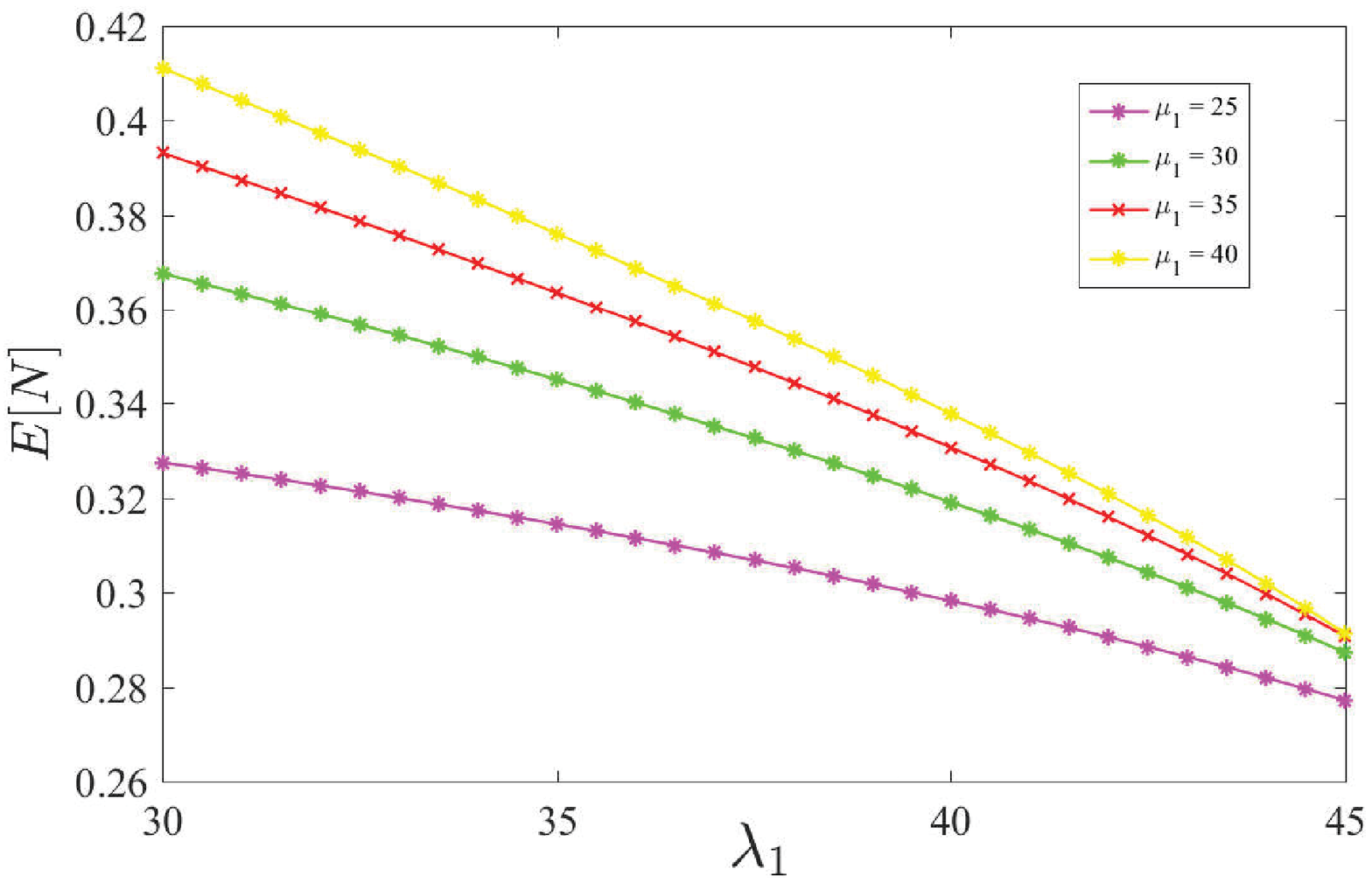} %
\includegraphics[width=7cm]{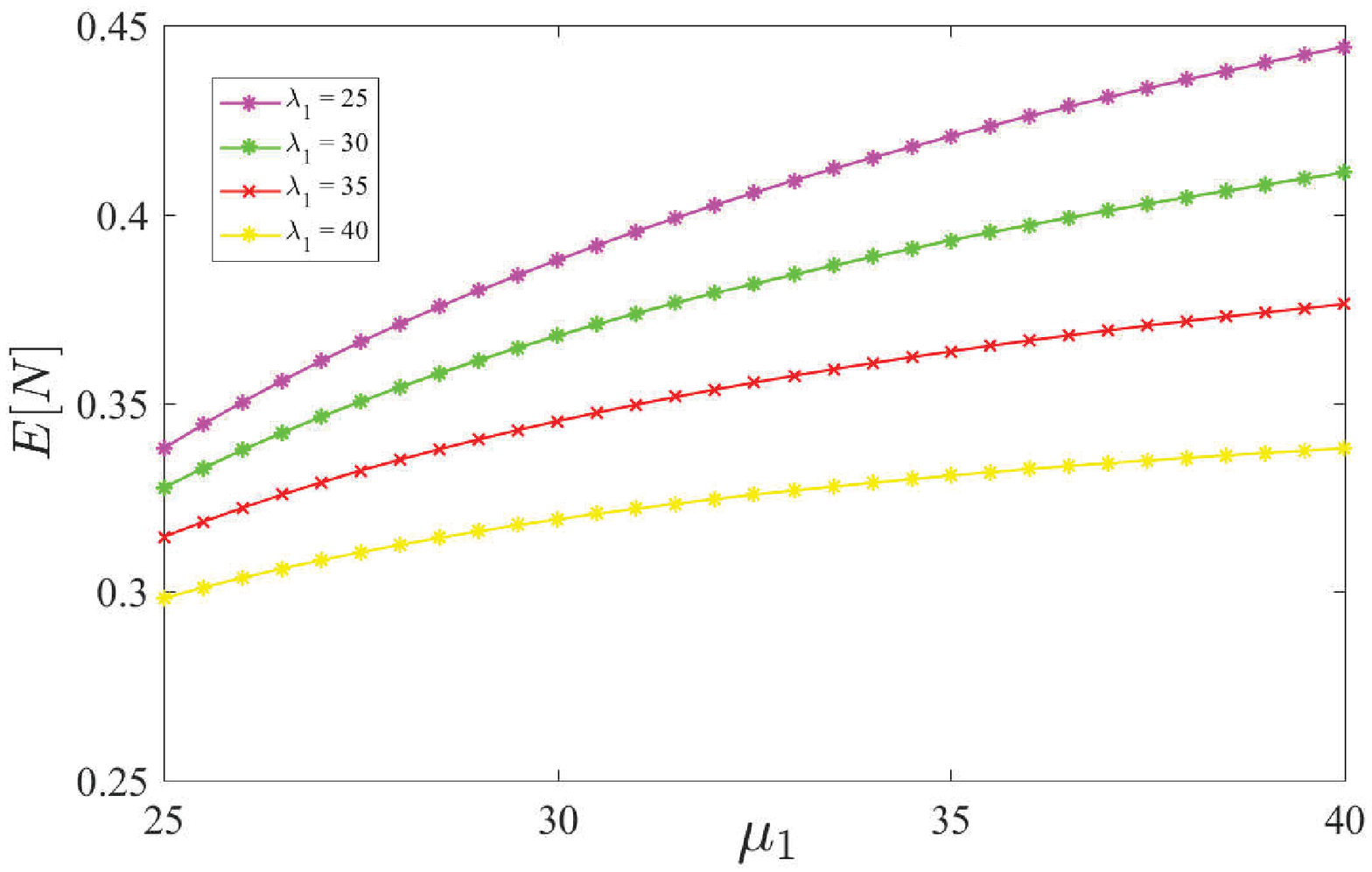} \newline
\caption{$E\left[ N\right] $ vs. $\protect\lambda_{1}$ and $\protect\mu_{1}$}
\label{figure:figure-6}
\end{figure}

\textbf{Example three: Analysis of the stationary weak-probability }$\mathbf{%
p}_{w}$

The left of Figure 7 shows how the stationary weak-probability $\mathbf{p}%
_{w}$ depends on $\lambda _{1}\in \left( 30,45\right) $ when $\mu _{1}=25$, $%
30$, $35$ and $40$, respectively. It is seen that $\mathbf{p}_{w}$ decreases
as $\lambda _{1}$ increases but it increases as $\mu _{1}$ increases. The
right of figure 7 shows how $\mathbf{p}_{w}$ depends on $\mu _{1}\in \left(
25,40\right) $ when $\lambda _{1}=25$, $30$, $35$ and $40$, respectively. It
is seen that $\mathbf{p}_{w}$ increases as $\mu _{1}$ increases but it
decreases as $\lambda _{1}$ increases. Note that the number of rented bikes
increases as $\lambda _{1}$ increases, thus $\sum_{k=-L}^{0}\pi _{k}e$
increases but $\sum_{k=K}^{K+L}\pi _{k}e$ decreases. On the other hand, the
number of returned bikes increases as $\mu _{1}$ increases. This indicates
that $\sum_{k=K}^{K+L}\pi _{k}e$ increases but $\sum_{k=-L}^{0}\pi _{k}e$
decreases. It is well understood that $\sum_{k=K}^{K+L}\pi _{k}e$ has more
impact on $\mathbf{p}_{w}$ than $\sum_{k=-L}^{0}\pi _{k}e$ under the present
parameter design.

\begin{figure}[ptb]
\setlength{\abovecaptionskip}{0.cm} \setlength{\belowcaptionskip}{-0.cm} %
\centering              \includegraphics[width=7cm]{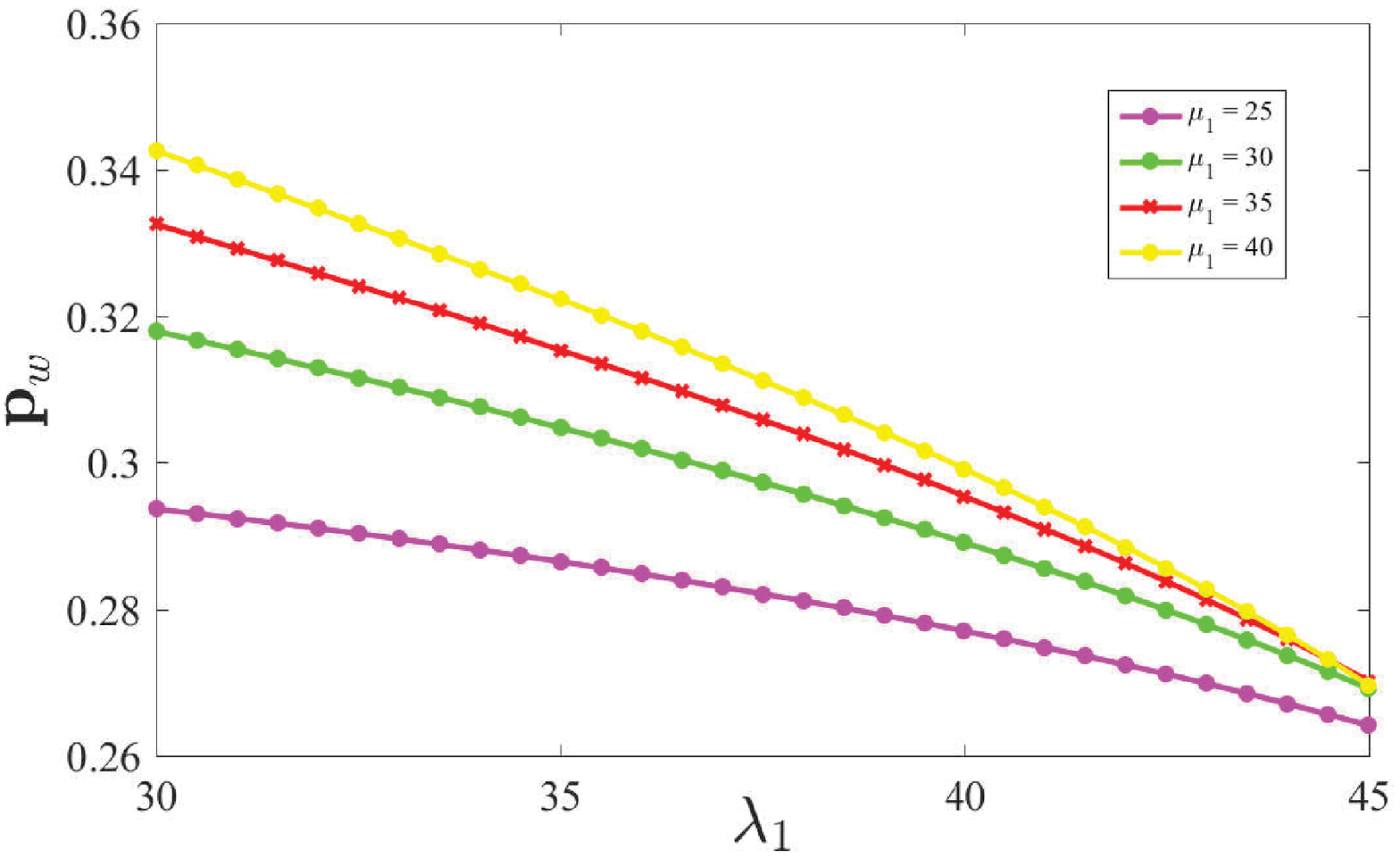} %
\includegraphics[width=7cm]{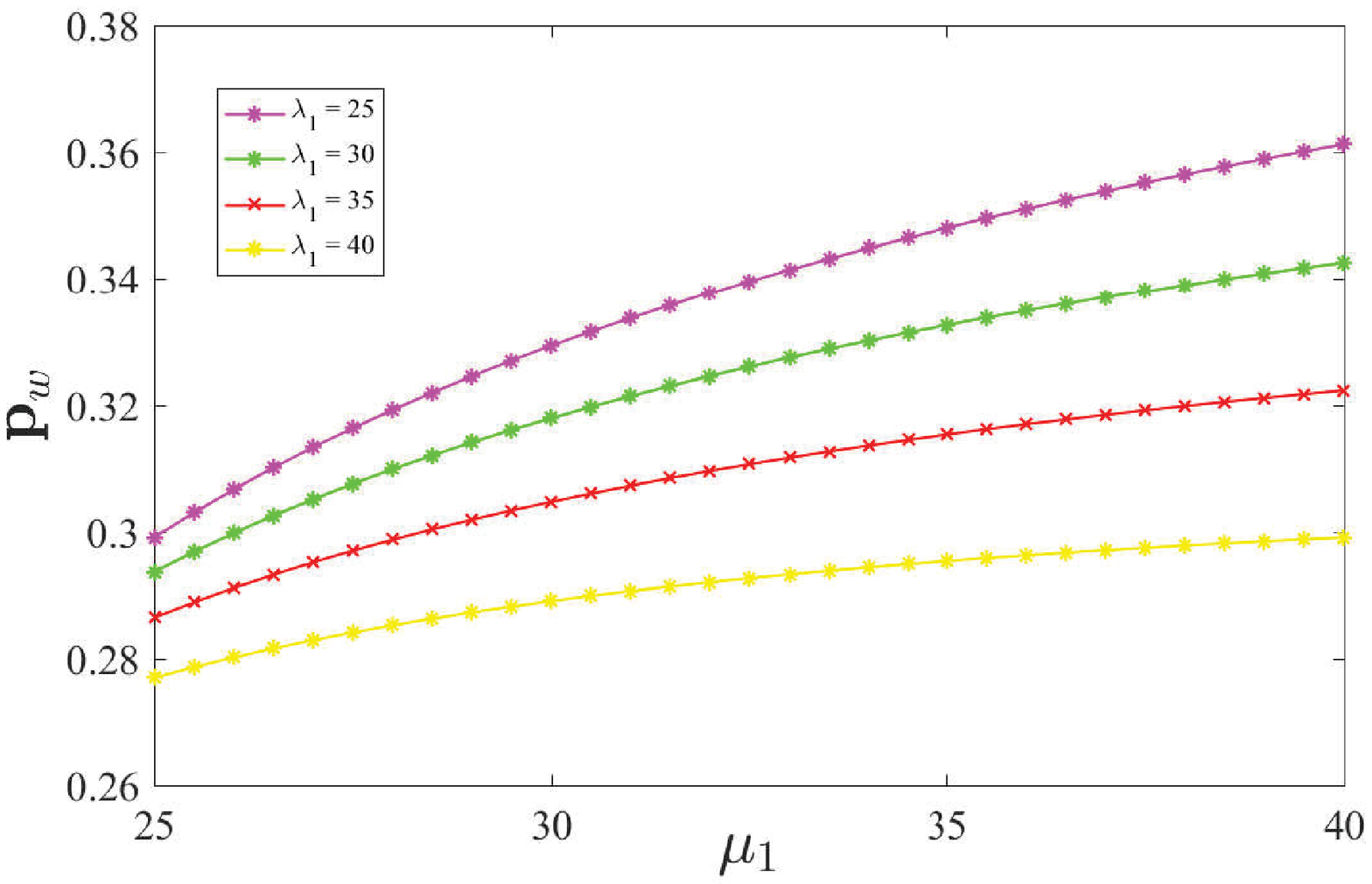} \newline
\caption{$\mathbf{p}_{w}$ vs. $\protect\lambda_{1}$ and $\protect\mu_{1}$}
\label{figure:figure-7}
\end{figure}

\textbf{Example four: Analysis of the stationary strong-probability }$%
\mathbf{p}_{s}$

The left of figure 8 shows how the stationary strong-probability\textbf{\ }$%
\mathbf{p}_{s}$ depends on $\lambda _{1}\in \left( 30,45\right) $ when $\mu
_{1}=25$, $30$, $35$ and $40$, respectively. It is seen that $\mathbf{p}_{s}$
decreases as $\lambda _{1}$ increases but it increases as $\mu _{1}$
increases. The right of figure 8 shows how $\mathbf{p}_{s}$ depends on $\mu
_{1}\in \left( 25,40\right) $ when $\lambda _{1}=25$, $30$, $35$ and $40$,
respectively. It is seen that $\mathbf{p}_{s}$ increases as $\mu _{1}$
increases but it decreases as $\lambda _{1}$ increases. Intuitively, the
number of rented bikes increases as $\lambda _{1}$ increases, thus $\pi
_{-L}e$ increases but $\pi _{K+L}e$ decreases. On the other hand, the number
of returned bikes increases as $\mu _{1}$ increases, hence $\pi _{K+L}e$
increases but $\pi _{-L}e$ decreases. This demonstrates that $\pi _{K+L}e$
has more impact on $\mathbf{p}_{s}$ than $\pi _{-L}e$.

\begin{figure}[tbp]
\setlength{\abovecaptionskip}{0.cm} \setlength{\belowcaptionskip}{-0.cm} %
\centering              \includegraphics[width=6.5cm]{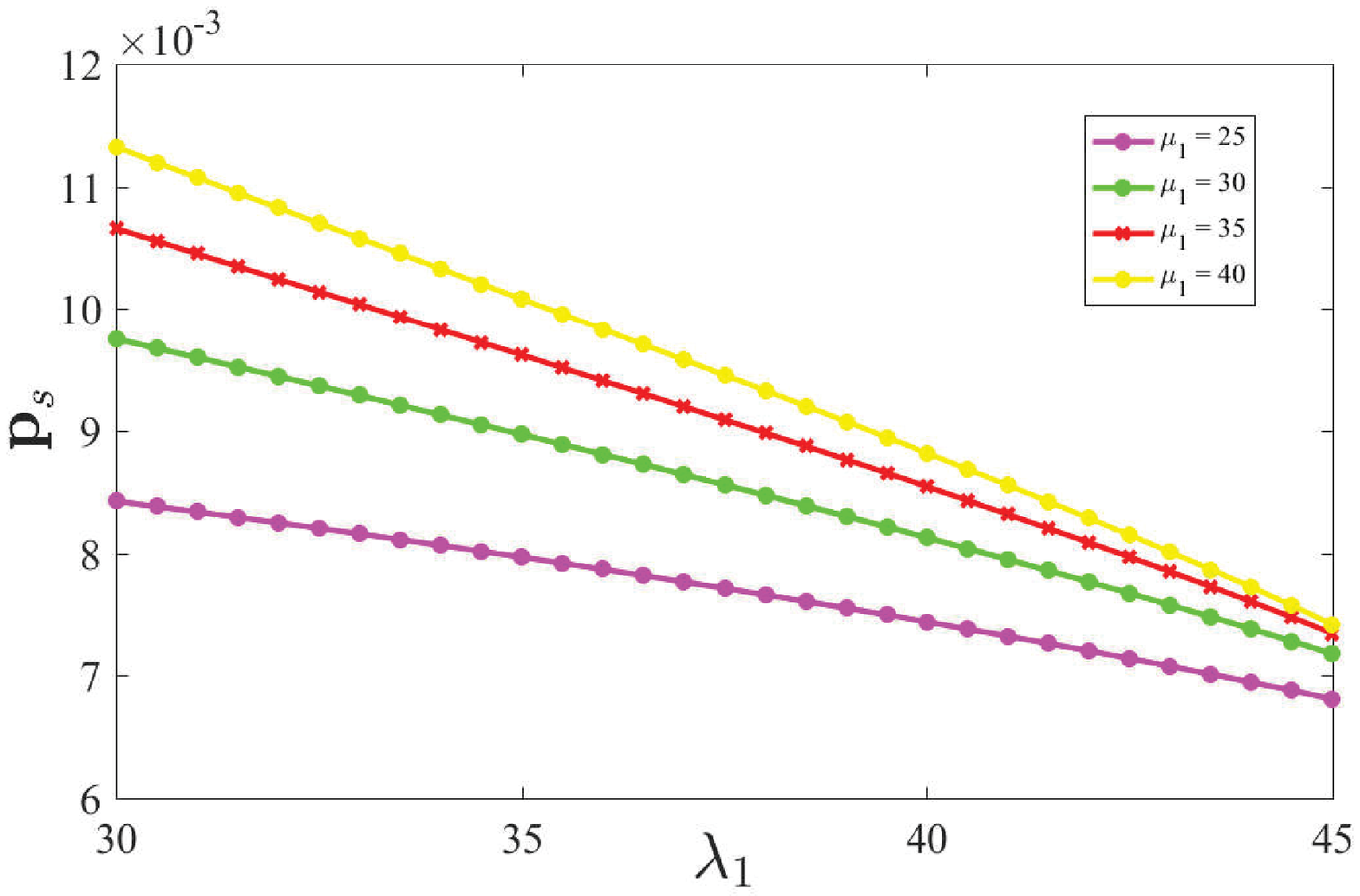} %
\includegraphics[width=7cm]{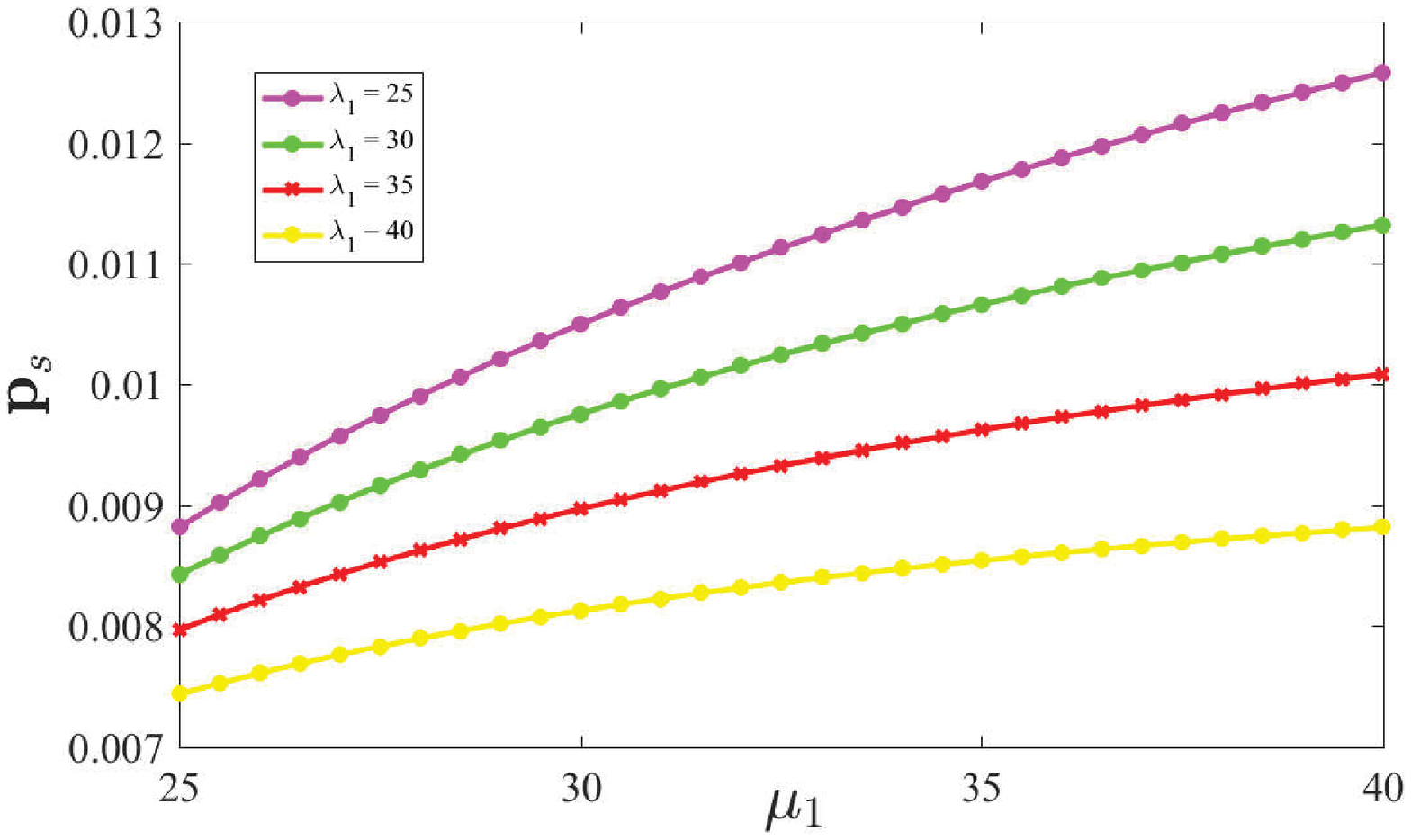} \newline
\caption{$\mathbf{p}_{s}$ vs. $\protect\lambda _{1}$ and $\protect\mu _{1}$}
\label{figure:figure-8}
\end{figure}

In the remainder of this section, we further observe some numerical impacts
of the user's finite waiting rooms on system performance through the
following two examples.

In Examples five and six, we take some common parameters as follows:
\begin{equation*}
K=20,C=5,m=2,\lambda _{2}=50,w=\left(
\begin{array}{cc}
-1 & 1 \\
1 & -1%
\end{array}%
\right) ;
\end{equation*}%
while the other parameters are conceretely chosen in each example for the
target of specific observation.

\textbf{Example five: Analysis of the stationary strong-probability }$%
\mathbf{p}_{s}$ \textbf{vs.} $L$

In the left of Figure 9, We take $\lambda _{1}=45$, $\mu _{1}=\mu _{2}=20$
and $\beta =0.75$. The left of Figure 9 shows how the stationary
strong-probability $\mathbf{p}_{s}$ depends on $L\in \left\{
0,1,2,3,4,5,6,7,8,9\right\} $ when $\alpha =0.10,0,25,0.50,0.75,0.90$,
respectively. It is seen that $\mathbf{p}_{s}$ decreases as $L$ increases
but it increases as $\alpha $ increases.

In the right of Figure 9, we take $\lambda _{1}=55$, $\mu _{1}=\mu _{2}=10$
and $\alpha =0.75$. The right of Figure 9 shows how $\mathbf{p}_{s}$ depends
on $L\in \left\{ 0,1,2,3,4,5,6,7,8,9\right\} $ when $\beta
=0.10,0,25,0.50,0.75,0.90$, respectively. It is seen that $\mathbf{p}_{s}$
decreases as $L$ increases but it increases as $\beta $ increases.

On the one hand, as $L$ increases, there are more and more waiting places
provided for customers to wait for either an available bike (rent) or a
vacant parking place (return), thus the probabilities of $\pi _{-L}e$ and $%
\pi _{K+L}e$ will decrease. This illustrates that $\mathbf{p}_{s}$
decreases. On the other hand, as $\alpha $ increases, more and more
customers would like to wait for an available bike, thus $\pi _{-L}e$
increases but $\pi _{K+L}e$ decreases. As $\beta $ increases, more and more
customers would like to wait for a vacant parking place, thus $\pi _{K+L}e$
increases but $\pi _{-L}e$ decreases.
\begin{figure}[tbp]
\setlength{\abovecaptionskip}{0.cm} \setlength{\belowcaptionskip}{-0.cm} %
\centering               \includegraphics[width=7cm]{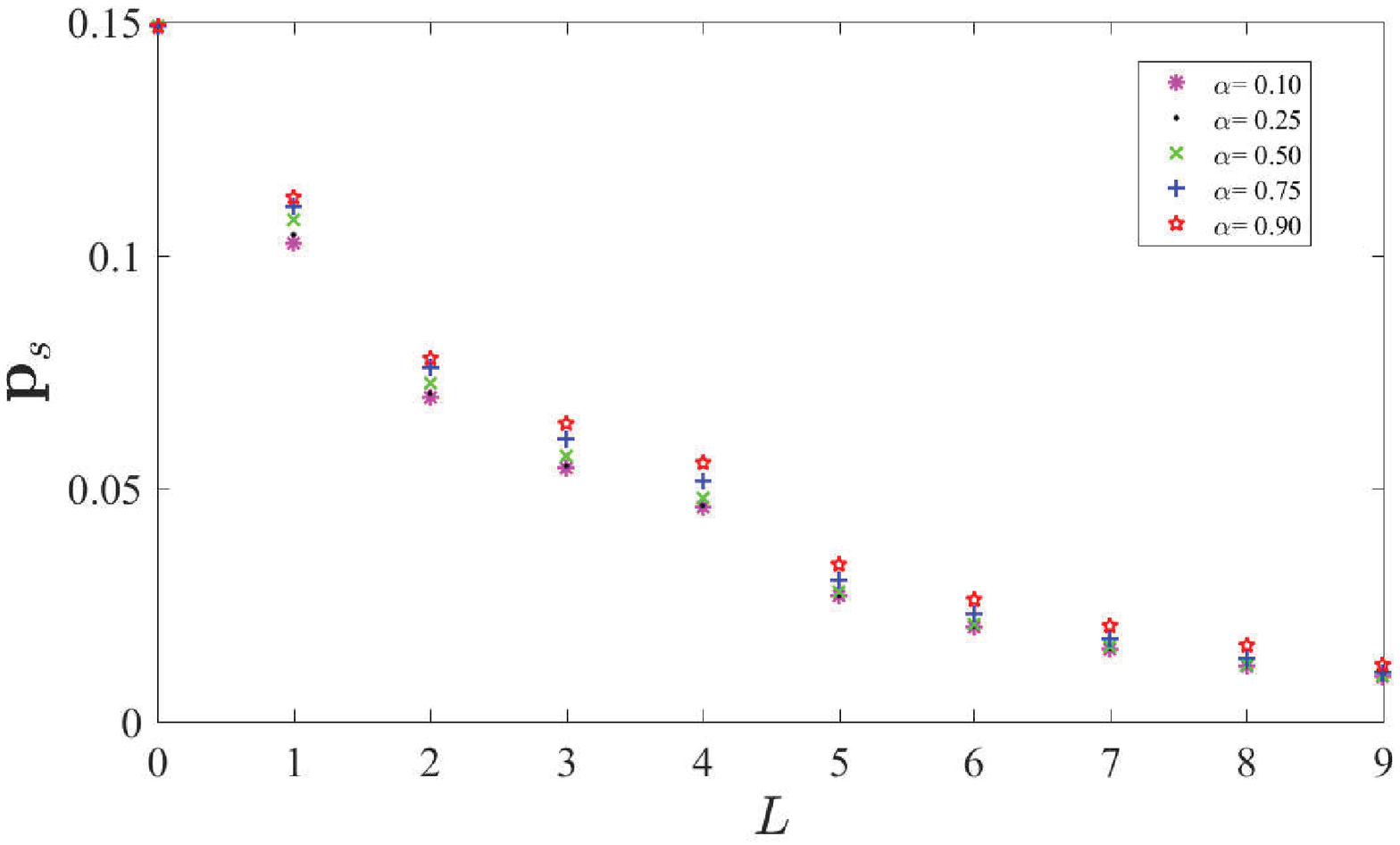} %
\includegraphics[width=7cm]{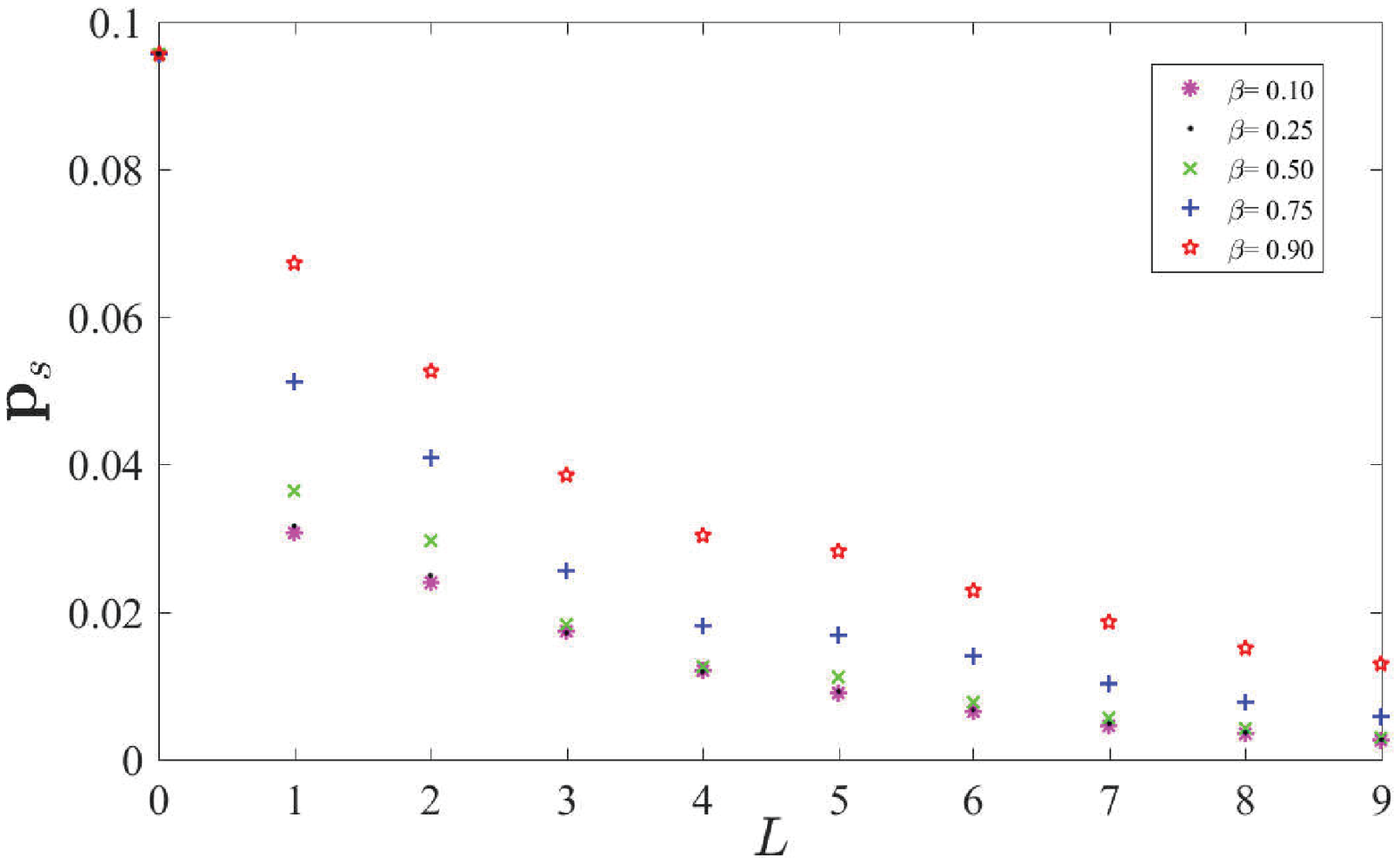} \newline
\caption{$\mathbf{p}_{s}$ vs. $L$}
\label{figure:figure-9}
\end{figure}

\textbf{Example six: Analysis of the ratio }$\upsilon $ \textbf{vs.} $L$

In the left of Figure 10, we take $\lambda _{1}=45$, $\mu _{1}=\mu _{2}=20$
and $\beta =0.75$. The left of Figure 10 shows how the ratio\textbf{\ }$%
\upsilon $ depends on $L\in \left\{ 0,1,2,3,4,5,6,7,8,9\right\} $ when $%
\alpha =0.10,0,25,0.50,0.75,0.90$, respectively. It is seen that the ratio $%
\upsilon $ decreases as $L$ increases but it increases as $\alpha $
increases.

In the right of Figure 10, we take $\lambda _{1}=55$, $\mu _{1}=\mu _{2}=10$
and $\alpha =0.75$. The right of Figure 10 shows how the ratio\textbf{\ }$%
\upsilon $ depends on $L\in \left\{ 0,1,2,3,4,5,6,7,8,9\right\} $ when $%
\beta =0.10,0,25,0.50,0.75,0.90$, respectively. It is seen that the ratio $%
\upsilon $ decreases as $L$ increases but it increases as $\beta $
increases.
\begin{figure}[tbp]
\setlength{\abovecaptionskip}{0.cm} \setlength{\belowcaptionskip}{-0.cm} %
\centering              \includegraphics[width=7cm]{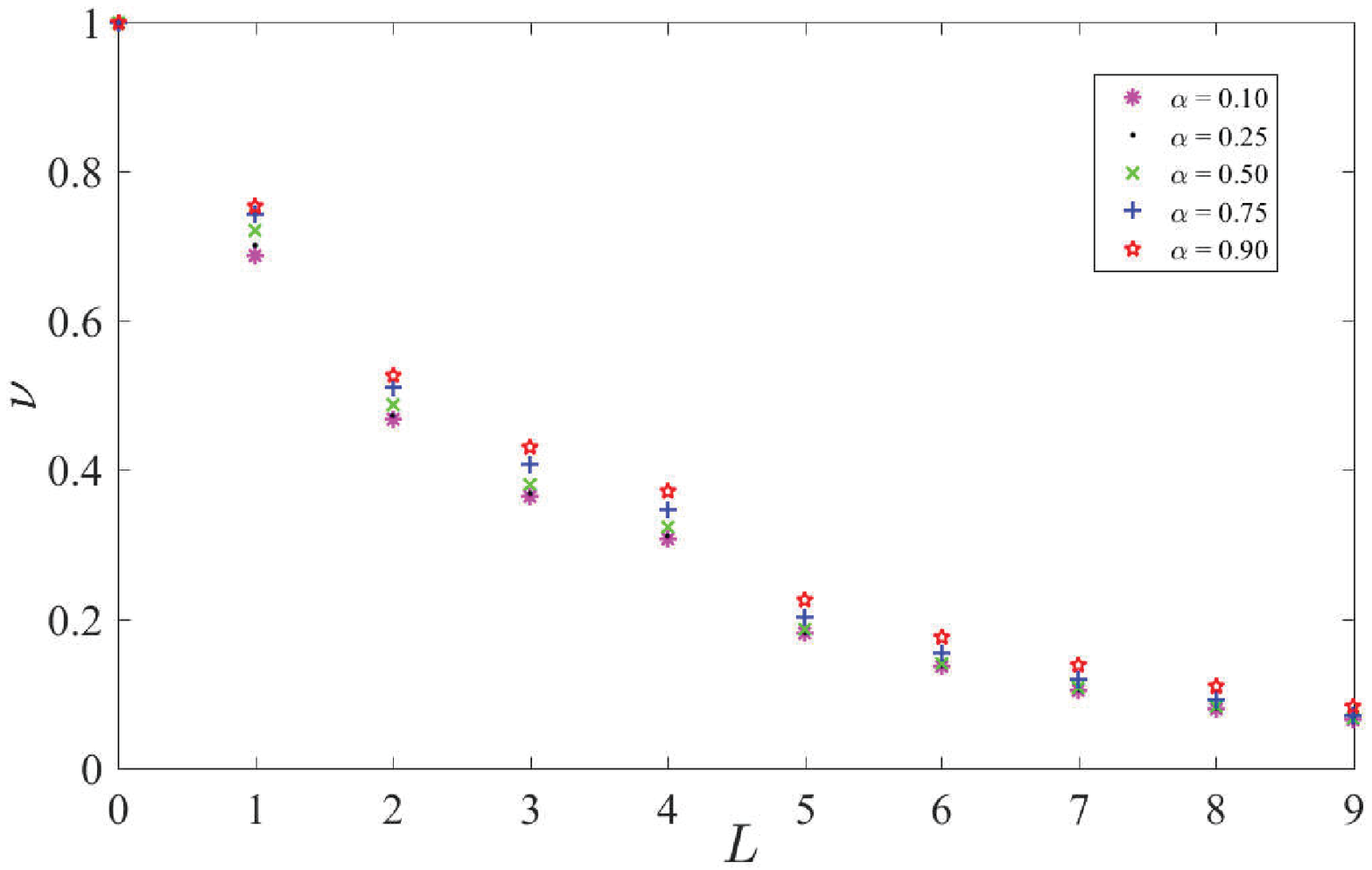} %
\includegraphics[width=7cm]{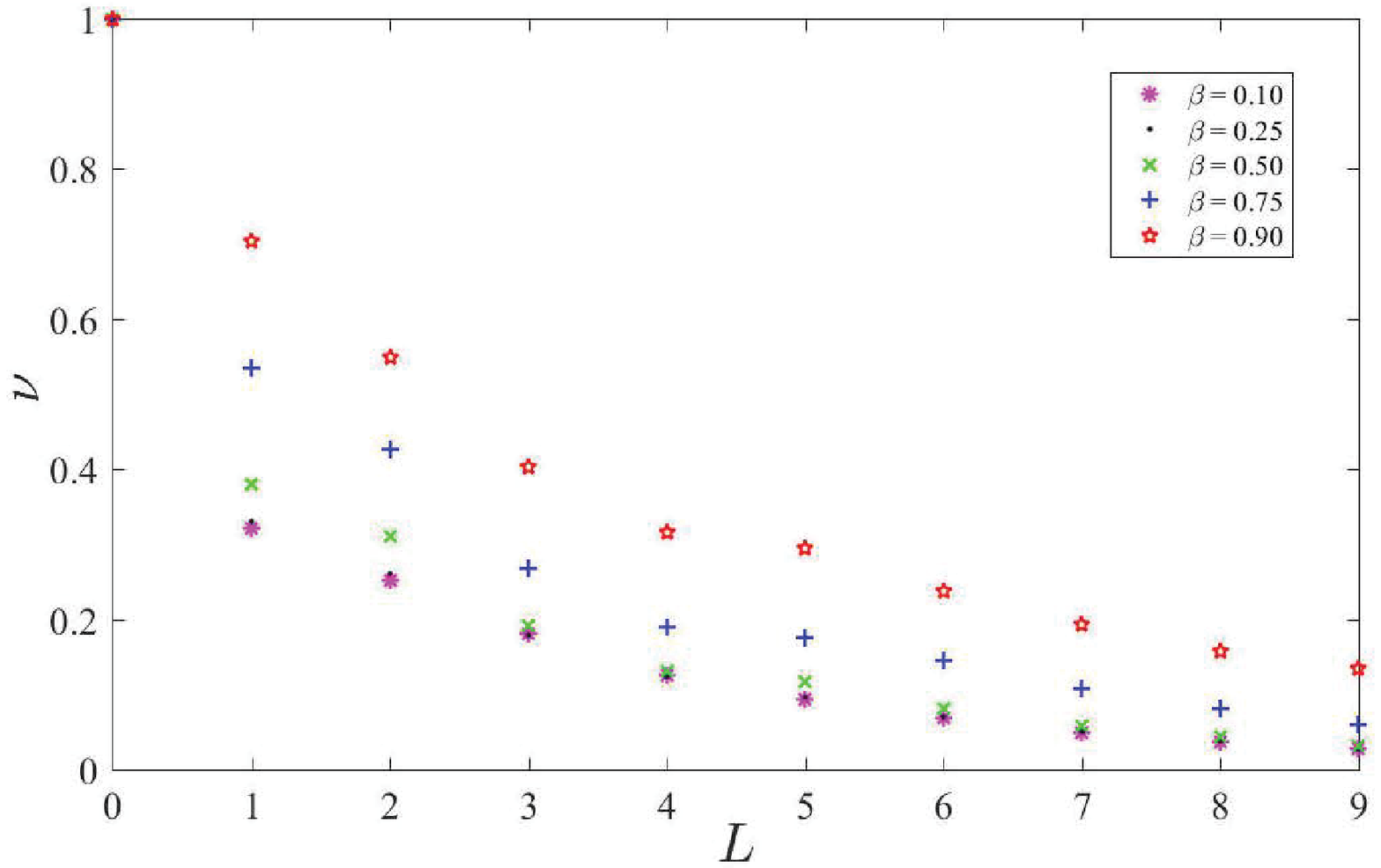} \newline
\caption{$\mathbf{p}_{s}$ vs. $\protect\upsilon $}
\label{figure:figure-10}
\end{figure}

\section{Concluding Remarks}

In this paper, we describe a large-scale bike sharing system under Markovian
environment, and develop a mean-field matrix-analytic method by combining
the mean-field theory with the time-inhomogeneous queues as well as the
nonlinear QBD processes. Furthermore, we apply the martingale limit theory
to prove the asymptotic independence (or propagation of chaos) of this bike
sharing system, and also study the limiting interchangeability as $%
N\rightarrow \infty $ and $t\rightarrow +\infty $. Based on this, we discuss
the fixed point by means of a nonlinear QBD process so that we can give
performance analysis of this bike sharing system. Notice that the mean-field
matrix-analytic method is effective and efficient for, such as, designing
reasonable architecture of a bike sharing system, finding a better path
scheduling, improving inventory management, redistributing the bikes among
stations or clusters in terms of truck scheduling, price optimization,
application of intelligent information technologies and so forth.

This paper provides a clear way for how to use the mean-field
matrix-analytic method to analyze performance measures of more general bike
sharing systems in practice through three key parts: (1) Setting up a
mean-field system of mean-field equations, (2) proving the asymptotic
independence, and (3) analyzing performance measures of this bike sharing
system by means of the fixed point. Therefore, the methodology and results
of this paper give some new highlight on understanding performance measures
and operations management of bike sharing systems. Along such a line, there
are a number of interesting directions for potential future research, for
example:

\begin{itemize}
\item Analyzing the fixed point for more general bike sharing systems in
practice, and provide effective algorithms to deal with the nonlinear QBD
processes;

\item studying non-exponential riding-bike times and non-Poisson customer
arrivals in bike sharing systems;

\item introducing some better operations management, such as, redistribution
of bikes by trucks, inventory management, applications of intelligent
information techniques; and

\item discussing large-scale bike sharing systems with different clusters
or/and under price optimization.
\end{itemize}

\section*{Acknowledgements}

Q.L. Li was supported by the National Natural Science Foundation of China
under grant No. 71471160 and No. 71671158, and Natural Science
Foundation of Hebei Province in China under grant No. G2017203277.

\end{document}